\documentclass[12pt]{amsart}
\usepackage{amsmath,amsfonts,amsthm,amscd}
\usepackage{graphicx}
\usepackage{esint}
\def\barint{\kern4pt
\raise3.4pt\hbox{\vrule height.8pt width5pt}%
\kern-9pt 
\int}
\newtheorem{theorem}{Theorem}

\newtheorem{proposition}[theorem]{Proposition}
\newtheorem{lemma}[theorem]{Lemma}

\newtheorem{corollary}[theorem]{Corollary}
\newtheorem{remark}[theorem]{Remark}

\renewcommand{\a }{\alpha }
\renewcommand{\b }{\beta }
\renewcommand{\d}{\delta }
\newcommand{\D }{\Delta }

\newcommand{\e }{\varepsilon }
\newcommand{\g }{\gamma}

\renewcommand{\l }{\lambda }
\renewcommand{\L }{\Lambda }

\newcommand{\n }{\nabla }
\newcommand{\var }{\varphi }

\newcommand{\pa }{\partial}

\newcommand{\ov}{\overline}

\newcommand{\be}{\begin{equation}}
\newcommand{\ee}{\end{equation}}
\newenvironment{pf}{\noindent{\sc Proof}.\enspace}{\rule{2mm}{2mm}\medskip}
\newenvironment{pfn}{\noindent{\sc Proof}}{\rule{2mm}{2mm}\medskip}
\newcommand{\R}{\mathbb{R}}

\newcommand{\vs}{\vskip.1in}

\numberwithin{equation}{section} \numberwithin{theorem}{section}
\begin{document}
\date{April 20, 2011}
\bibliographystyle{amsalpha}

\title[non-uniqueness]{Non-uniqueness results for critical metrics of regularized determinants in four dimensions}

\author{Matthew Gursky}
\address{Department of Mathematics\\
         University of Notre Dame\\
         255 Hurley Hall \\
         Notre Dame, IN 46556}
\email{mgursky@nd.edu}

\author{Andrea Malchiodi}
\address{Sector of Mathematical Analysis \\
         SISSA  \\
         Via Bonomea 265  \\
         34136 Trieste \\
         ITALY}
\email{malchiod@sissa.it}

\thanks{First author supported in part by NSF grant DMS-0800084.  Second
author supported by the FIRB project {\em Analysis and Beyond} from MIUR. Both authors
are grateful to IAS in Princeton and to the De Giorgi Center in Pisa for the kind
hospitality during the preparation of this work. The authors would like to thank Prof.
Arthur Lim, who produced the graphics of the invariant disc in Figures 3, 5 and 7.}

\maketitle

\begin{abstract}    The regularized determinant of the Paneitz operator arises in
quantum gravity (see \cite{Connesbook}, IV.4.$\gamma$).
An explicit formula for the relative determinant of two conformally related metrics
was computed by Branson in \cite{BransonCMP}.  A similar formula holds for Cheeger's half-torsion,
which plays a role in self-dual field theory (see \cite{Juhlbook}), and is defined in terms of regularized
determinants of the Hodge laplacian on $p$-forms ($p < n/2$).  In this
article we show that the corresponding actions are unbounded (above and below) on any
conformal four-manifold.  We also show that the
conformal class of the round sphere admits a second solution which is not given by
the pull-back of the round metric by a conformal map,
thus violating uniqueness up to gauge equivalence.  These results differ from the
properties of the determinant of the conformal Laplacian
established in \cite{ChangYangAnnals}, \cite{BCY}, and \cite{Gursky2}.

We also study entire solutions of the Euler-Lagrange equation of $\log \det P$ and the half-torsion $\tau_h$ on
$\mathbb{R}^4 \setminus \{ 0 \}$, and show the existence of two
families of periodic solutions.  One of these families includes {\em Delaunay}-type
solutions.
\end{abstract}

\section{Introduction}

Let $(M^n,g)$ be a closed Riemannian manifold.  Let $\Delta =
\Delta_g$ denote the Laplace-Beltrami operator, and label the
eigenvalues of $(-\Delta_g)$ by
\begin{align*}
0 = \lambda_0 < \lambda_1 \leq \lambda_2 \leq \dots
\end{align*}
counting multiplicities.
The {\em spectral zeta
function} of $(M^n,g)$ is
\begin{align} \label{zetadef}
\zeta(s) = \sum_{j=1}^{\infty} \lambda_j^{-s}. 
\end{align}
By Weyl's asymptotic law,
\begin{align*}
\lambda_j \sim j^{2/n}, \ j \rightarrow \infty.
\end{align*}
Consequently, (\ref{zetadef}) defines an analytic function for
$\mbox{Re}(s) > n/2$.

Note that formally--that is, if we were to take the definition
in (\ref{zetadef}) literally--then
\begin{align} \label{formdet}
\zeta^{\prime}(0) = - \sum_{j = 1}^{\infty} \log \lambda_j = - \log \det (-\Delta_g),
\end{align}
although of course the series (\ref{zetadef}) does not
define an analytic function near $s = 0$.  However, one can meromorphically extend
so that $\zeta$ becomes regular at $s = 0$ (see
\cite{RS}), and in view of (\ref{formdet}) define the regularized determinant by
\begin{align} \label{logdetdef}
\det(-\Delta_g) = e^{- \zeta^{\prime}(0)}.
\end{align}

Since the determinant is obviously a global invariant, it is all the more remarkable that
Polyakov was able to write a local formula (appearing as a partition function in string
theory) for the ratio of the
determinants for two conformal metrics on a closed surface (see \cite{Polyakov}).
Suppose $\hat{g} = e^{2w}g$, then
\begin{align} \label{Polyform}
\log \frac{\det (-\Delta_{\hat{g}})}{\det (-\Delta_g)} =
-\frac{1}{12\pi} \int_{\Sigma} (|\nabla w|^2 + 2Kw)\ dA,
\end{align}
where $K = K_{g}$ is the Gauss curvature of $g$.

The formula (\ref{Polyform}) defines an action on the space of unit
volume conformal metrics $[g]_1 = \{ e^{2w}g\ |\ Vol(e^{2w}g) = \int
e^{2w}\ dA = 1 \}$.  Critical points of this action are precisely
those metrics of constant Gauss curvature; to see this one appeals
to the {\em Gauss curvature equation}
\begin{align} \label{GCE}
\Delta w + K_{\hat{g}} e^{2w} = K,
\end{align}
and computes a first variation of (\ref{Polyform}). In a series of
papers \cite{OPS1}, \cite{OPS2}, Osgood-Phillips-Sarnak studied the
existence of extremals for this functional, and the beautiful connection to various sharp
Moser-Trudinger-Sobolev inequalities.

\subsection{Four dimensions}

In deriving (\ref{Polyform}) Polyakov exploited a crucial property
of the Laplacian in two-dimensions, namely, its conformal
covariance: if $\hat{g} = e^{2w}g$, then
\begin{align*}
\Delta_{\hat{g}} = e^{-2w}\Delta_g.
\end{align*}
In general, we say that the metric-dependent
differential operator $A = A_g$ is {\em conformally covariant of
bi-degree $(a,b)$} if $\hat{g} = e^{2w}g$ implies
\begin{align} \label{Ccovar}
A_{\hat{g}}\psi = e^{-bw} A_g (e^{aw}\psi)
\end{align}
for each smooth section $\psi$ of some vector bundle $\mathbb{E}$.  Examples of such operators include the
conformal Laplacian
\begin{align} \label{Ldef}
L = -\Delta + \frac{(n-2)}{4(n-1)}R,
\end{align}
where $R$ is the scalar curvature, with $a = \frac{n-2}{2}$ and $b = \frac{n+2}{2}$, and the
four-dimensional Paneitz operator
\begin{align} \label{Pdef}
P = (-\Delta)^2 + \delta \left( \frac{2}{3}Rg - 2 Ric \right) \circ \nabla,
\end{align}
with $a = 0$ and $b = 4$.  Indeed, the Paneitz operator is from many points
of view the
natural generalization of the Laplace-Beltrami operator to
four-manifolds, and in analogy to the Gauss curvature equation we
have the prescribed $Q$-curvature equation
\begin{align} \label{QCE}
Pw + 2 Q = 2Q_{\hat{g}}e^{4w},
\end{align}
where $Q$ is the $Q$-curvature:
\begin{align} \label{Qdef}
Q = \frac{1}{12}( -\Delta R + R^2 - 3|Ric|^2).
\end{align}

In \cite{BO}, Branson-$\O$rsted were able to generalize Polyakov's
technique to conformally covariant operators $A$ defined on a
four-manifold $M^4$. The resulting formula, while somewhat
complicated, is geometrically quite natural.  The first thing to
note is that it is always a linear combination of three universal terms appearing
in the determinant formula, with different linear combinations depending on the
choice of operator $A$.  Therefore, the formula is typically
expressed as
\begin{align} \label{FAdef}
F_{A}[w] = \log \frac{\det A_{\hat{g}}}{\det A_g} = \gamma_1(A) I[w]
+ \gamma_2(A) II[w] + \gamma_3(A) III[w],
\end{align}
where $(\gamma_1,\gamma_2,\gamma_3)$ is a triple of real numbers,
and $I,II,III$ are the three {\em sub-functionals}.  For example, if $A
= L$, the conformal Laplacian, then
\begin{align} \label{gammaL}
\gamma_1(L) = 1, \qquad \gamma_2(L)= -4, \qquad \gamma_3(L) = -2/3.
\end{align}

In general, if $A$ has a non-trivial kernel, then one needs to
modify the definition of the zeta function (since $0$ is an eigenvalue);
this results in some additional terms in the formula for $F_A$, see \cite{BO}.

Before giving the precise formulas for these functionals, it may
shed some light if we first describe their geometric content:
\begin{align*} 
\hat{g} = e^{2w}g\ \mbox{is a critical point of }I\
\Longleftrightarrow |W_{\hat{g}}|^2 = const.,
\end{align*}

\begin{align*} 
\hat{g} = e^{2w}g\ \mbox{is a critical point of }II\
\Longleftrightarrow Q_{\hat{g}} = const.,
\end{align*}

\begin{align*} 
\hat{g} = e^{2w}g\ \mbox{is a critical point of }III\
\Longleftrightarrow \Delta_{\hat{g}} R_{\hat{g}} = 0,
\end{align*}
where $W$ is the Weyl curvature tensor. Thus, each functional
corresponds to a natural curvature condition in four dimensions.
The functionals $II$ and $III$ are of particular interest as they
correspond to, respectively, the constant $Q$-curvature problem and
the Yamabe problem. \vskip.2in

\subsection{The formulas} The precise formulas\footnote{In fact, Branson-$\O$rsted considered a
scale-invariant version of the regularized determinant;
hence each functional above is invariant under $w \mapsto w + c$.} for $I$, $II$, and $III$ are
\begin{align} \label{Idef}
I[w] = 4 \int w| W |^2\ dv - \big( \int |W|^2\ dv \big) \log \fint
e^{4w}\ dv,
\end{align}
\begin{align} \label{IIdef}
II[w] = \int w Pw\ dv - \big( \int Q\ dv\big) \log
\fint e^{4(w - \overline{w})}\ dv,
\end{align}
\begin{align} \label{IIIdef}
III[w] = 12 \int (\Delta w + |\nabla w|^2 )^2\ dv - 4 \int ( w
\Delta R + R|\nabla w|^2 )\ dv.
\end{align}

In order to write down the Euler-Lagrange equation for $F_A$, we
define the following conformal invariant:
\begin{align} \label{kappadef}
\kappa_A = -\gamma_1 \int |W|^2\ dv - \gamma_2 \int Q\ dv.
\end{align}
Then the E-L equation is
\begin{align} \label{EL1} \begin{split}
\mu e^{4w} &= (\frac{1}{2}\gamma_2 + 6\gamma_3)\Delta^2 w +
6\gamma_3 \Delta |\nabla w|^2 - 12\gamma_3 \nabla^i \big[ (\Delta w
+ |\nabla w|^2)\nabla_i w\big] \\
& \hskip.25in + \gamma_2 R_{ij}\nabla_i \nabla_j w + (2\gamma_3 -
\frac{1}{3}\gamma_2)R\Delta w + (2\gamma_3 + \frac{1}{6}\gamma_2)
\langle \nabla R,\nabla w \rangle \\
& \hskip.3in + (\gamma_1 |W|^2 +\gamma_2 Q - \gamma_3 \Delta R),
\end{split}
\end{align}
where
\begin{align} \label{mudef}
\mu = -\frac{\kappa_A}{\int e^{4w}}.
\end{align}
Note the equations are in general fourth order, unless $\frac{1}{2}\gamma_2 + 6\gamma_3 = 0$.  In this case
the equation is second order but fully nonlinear; it is precisely the $\sigma_2$-curvature condition (see \cite{JeffSurvey}).

Geometrically, (\ref{EL1}) means the following:  denote the {\em $U$-curvature} of $g$
\begin{align} \label{Udef}
U = U(g) = \gamma_1 |W|^2 +\gamma_2 Q - \gamma_3 \Delta R.
\end{align}
If $w$ satisfies (\ref{EL1}), then the conformal metric $g_A =
e^{2w}g$ satisfies
\begin{align} \label{EL2}
U(g_A) \equiv \mu.
\end{align} \vskip.2in

\subsection{Some general existence results}

The first existence results for extremals of the functional
determinant in four dimensions were proven by Chang-Yang
\cite{ChangYangAnnals}.

\begin{theorem} \label{thm1} {\em (Chang-Yang, \cite{ChangYangAnnals})}
Assume:

\vskip.1in \noindent $(i)$ $\gamma_2 < 0$ and $\gamma_3 < 0$,

\vskip.1in \noindent $(ii)$ $\kappa_A < (-\gamma_2)8 \pi^2$.

\vskip.1in  Then $\sup_{w \in W^{2,2}} F_A[w]$ is attained by some
$w \in W^{2,2}$.
\end{theorem}
\vskip.1in

For example, taking $A = L$ the conformal
Laplacian, then an extremal exists for $F_L$ provided
\begin{align} \label{subcrit}
\kappa_{L}(M^4,g) = -\int |W|^2\ dv + 4 \int Q\ dV < 32\pi^2.
\end{align}
This condition is related to the best constant in the Moser-Trudinger
inequality of Adams \cite{Adams}, and eliminates the possibility of
bubbling (note for the round sphere, $\kappa_L = 32\pi^2$). Regularity
of extremals was proved by the first author in joint work with Chang-Yang
\cite{CGYAJM}; later Uhlenbeck-Viaclovsky proved a more general regularity result for any
critical point of (\ref{EL1}) (see \cite{UVMRL}).

The first author established that the
condition (\ref{subcrit}) is always satisfied by a $4$-manifold of
non-negative scalar curvature, unless it is conformally equivalent
to the round sphere \cite{Gursky1}.  In this case, Branson-Chang-Yang proved that the round metric, and its
orbit under the conformal group, maximizes $F_{L}$ \cite{BCY}.
Later, the first author proved that the round metric (modulo the conformal
group) is the {\em unique} critical point \cite{Gursky2}.  Thus the existence theory for $F_L = \log \det L$, at least
for $4$-manifolds of positive scalar curvature, is complete, and we have uniqueness (modulo the conformal group)
on the sphere.

In general situations not much is known about existence of critical points. In \cite{DM}
the functional $II$ is studied in generic situations, and saddle points solutions
are found using a global variational scheme.

\subsection{Determinant of the Paneitz operator and Cheeger's half-torsion}

In this paper we are interested in regularized determinants for which condition $(i)$ of Theorem \ref{thm1}
fails; i.e., the coefficients $\gamma_2$ and $\gamma_3$ have
different signs.  The corresponding functionals are therefore non-convex combinations of terms with different homogeneities,
and their variational properties quite difficult to analyze.  This arises in two cases of interest in mathematical physics:
the determinant of the Paneitz operator, and {\em Cheeger's half-torsion}.

In his book {\em Noncommutative Geometry}, Alain Connes devoted a
section to the discussion of the determinant of the Paneitz operator
(see \cite{Connesbook}, Chapt. IV.4.$\gamma$), ending with the
remark that "...the gravity theory induced from the above scalar
field theory in dimension $4$ should be of great interest..." In \cite{BransonCMP}, Branson
calculated the coefficients of $F_P$ and found $(\gamma_1,\gamma_2,\gamma_3) = (-1/4,-14,8/3)$.

For even-dimensional manifolds the {\em half-torsion} is defined by
\begin{align} \label{ht}
\tau_h = \frac{ \displaystyle (\det(-\Delta_0))^n (\det (-\Delta_2))^{n-4} \dots}{\displaystyle (\det (-\Delta_1))^{n-2} (\det (-\Delta_3))^{n-6}\cdots},
\end{align}
where $\Delta_p$ denotes the Hodge laplacian on $p$-forms.  Notice that this only involves $p$ for $p < n/2$; in particular
in four dimensions we have
\begin{align} \label{tau4}
\tau_h = \frac{ \displaystyle (\det(-\Delta_0)^4 }{\displaystyle (\det(-\Delta_1))^2}.
\end{align}
The half-torsion plays a role in self-dual field theory, for which the dimensions of physical interest
are $n = 4\ell + 2$. 
Witten's novel approach to studying self-dual field theory involved using Chern-Simons
theory in $4 \ell+3$-dimensions (see \cite{Witten97}). Cheeger's half-torsion appears when
computing the metric dependence of the partition function, similar
to Polyakov's formula (\cite{BelovMoore06}, \cite{Monnier10}).  Note that although the Hodge laplacian in general
does does not
satisfy (\ref{Ccovar}), the ratio in (\ref{ht}) has the requisite conformal properties for deriving a Polyakov-type formula
 (see \cite{Juhlbook}, Section 6.15).  The coefficients for the corresponding functional are $(\gamma_1,\gamma_2,\gamma_3) = (-13,-248,116/3)$.

In this paper we consider these functionals in the case of the round $4$-sphere.  For the determinant of
the Paneitz operator we have
\begin{align} \label{FPS4} \begin{split}
F_P[w] &= \int \Big[ 18(\Delta w)^2 + 64 |\nabla w|^2 \Delta w + 32
|\nabla w|^4 - 60 |\nabla w|^2 \Big]\ dv \\ & \hskip.5in + 112\pi^2
\log \Big( \fint e^{4(w - \overline{w})}\ dv \Big).
\end{split}
\end{align}
Notice the cross term $\Delta w |\nabla w|^2$, and the fact that the
coefficient of $64$ is too large to allow this term to be absorbed
into the other (positive) terms.
Similarly, for the half-torsion we have
\begin{align} \label{FtS4} \begin{split}
F_{\tau}[w] &= \int \Big[ 216 (\Delta w)^2 + 928 |\nabla w|^2 \Delta w + 464
|\nabla w|^4 - 2352 |\nabla w|^2 \Big]\ dv \\ & \hskip.5in + 1984 \pi^2
\log \Big( \fint e^{4(w - \overline{w})}\ dv \Big).
\end{split}
\end{align}
Again, the exponential term has a 'good' sign, while the cross term can dominate the
other leading terms.
Compare these with the formula for
the determinant of $L$:
\begin{align} \label{FLS4} \begin{split}
F_L[w] &= \int \Big[ -12(\Delta w)^2 - 16 |\nabla w|^2 \Delta w -8
|\nabla w|^4 + 24 |\nabla w|^2 \Big]\ dv \\ & \hskip.5in + 32 \pi^2
\log \Big( \fint e^{4(w - \overline{w})}\ dv \Big).
\end{split}
\end{align}
In this case, the cross term can be absorbed into the other (negative) terms,
so the difficulty in proving the boundedness of a maximizing sequence is understanding
the interaction of the derivative terms with the exponential
term (this is precisely where the sharp inequality of Adams becomes crucial).

The Euler-Lagrange equation associated to (\ref{FPS4}) is
\begin{align} \label{ELS4} \begin{split}
 -42 e^{4w} &= 9 \Delta^2 w + 32 |\nabla^2 w|^2 - 32 (\Delta w)^2 - 32
\Delta u\ |\nabla u|^2 -32 \langle \nabla w, \nabla |\nabla w|^2 \rangle \\
&\hskip.25in + 78 \Delta u + 96 |\nabla w|^2 - 42.
\end{split}
\end{align}
Therefore, $w = 0$ (the round metric) is a critical point. In \cite{BransonCMP} Branson
calculated the second variation at $w=0$ and showed that it was a local minimum
(modulo deformations generated by the conformal group and rescalings). A similar calculation
shows that $w=0$ is a local minimum of $F_{\tau}$.   However, globally $F_{P}$ and $F_{\tau}$ are
{\em never} bounded from below:

\begin{theorem} \label{Thm1}
If $(M^4,g)$ is a closed four-manifold, then
\begin{align*}
\inf_{w \in W^{2,2}} F_{P}[w] &= -\infty, \\
\sup_{w \in W^{2,2}} F_{P}[w] &= +\infty,
\end{align*}
and the same holds for $F_{\tau}$.
\end{theorem}

While this result rules out using the direct approach for
finding critical points of $F_P$ or $F_{\tau}$, Branson's calculation suggests the possibility
of locating a second solution by looking for saddle points, for example by using
the Mountain Pass Theorem.  Of course, the conformal invariance of the functionals
implies that the Palais-Smale condition does not hold, so we need to somehow
mod out by the action of the conformal group, for example by imposing a symmetry condition.

The main result of this paper is the existence of a second (non-equivalent) critical point for $F_P$ and $F_{\tau}$ in the conformal class
of the round metric:

\begin{theorem} \label{Main}  Let
\begin{align*}
\mathbb{S}^4 = \{ (x_1,\dots,x_5) \in \mathbb{R}^5\ : x_1^2 + \cdots + x_5^2 = 1 \}
\end{align*}
be the $4$-sphere, and $g_0$ the round metric it inherits as a submanifold of $\mathbb{R}^5$.
Then there is critical point $u_P \in C^{\infty}(\mathbb{S}^4)$ of
$F_P$ such that  \vskip.1in

\noindent $(i)$  $u_P$ is rotationally symmetric and even:
\begin{align*}
u_P = u_P(x_5),\ u_P(x_5) = u_P(-x_5).
\end{align*}
\vskip.1in

\noindent $(ii)$ The metric $g = e^{2u_P}g_0$ is not conformally equivalent to $g_0$; i.e, there is
no conformal map $\varphi : \mathbb{S}^4 \rightarrow \mathbb{S}^4$ with $\varphi^{\ast}g = g_0$.
\vskip.1in

Moreover, $F_{\tau}$ admits a second solution $u_{\tau} = u_{\tau}(x_5)$ which is rotationally symmetric, even, but not
conformally equivalent to the round metric.
\end{theorem}
\vskip.1in

\noindent {\bf Remarks.} \vskip.1in

\noindent {\bf 1.}
In both cases, rotational symmetry reduces the Euler equation to an ODE.  Since the cylinder is conformal to the
sphere minus two points, we look for solutions on $\mathbb{R}^{+}$ with the appropriate
asymptotic behavior at infinity; see Section \ref{Y}.  \vskip.1in

\noindent {\bf 2.} The claim $(ii)$ of non-equivalence is actually immediate from the symmetry condition in $(i)$, since
evenness is not preserved by the action of the conformal group. \vskip.2in

In principle, one could exploit the variational structure of the problem and try
to apply standard variational methods like the Mountain Pass theorem. However it
seems difficult (even restricting to symmetric functions) to derive a-priori estimates
in $W^{2,2}$ on solutions or on Palais-Smale sequences, namely sequences of functions
satisfying
$$
F_P[u_k] \to c \in \R, \qquad \qquad  F'_P[u_k] \to 0,
$$
and similarly for $F_\tau$. For Yamabe-type problems, see e.g. \cite{BN}, to tackle the loss of
compactness one can first use energy bounds and classification of blow-up profiles,
which are lacking at the moment in our case.  \vskip.2in

It strikes us as somewhat remarkable that the sphere should admit a second distinct solution. Of course,
there is an abundance of examples in the literature in which the variational structure of an equation is
exploited to prove multiplicity results; but we are unaware of any geometric variational problems for
which constant curvature (mean, scalar, $Q$) does not characterize the sphere up to equivalence.
\vskip.2in

%



\subsection{Entire solutions.} A related question is the existence of solutions to the Euler equation for
$F_P$ or $F_{\tau}$ on Euclidean space.  For $F_P$ the equation is
\begin{align} \label{EulerR4}
c e^{4w} = 9 \Delta^2 w + 32 |\nabla^2 w|^2 - 32 (\Delta w)^2 - 32 \Delta u\ |\nabla u|^2 -32 \langle \nabla w, \nabla |\nabla w|^2 \rangle,
\end{align}
where $c$ is a constant (compare with (\ref{EL1})).  For $F_{\tau}$ we have
\begin{align} \label{EulerR4tau}
c' e^{4w} = 108 \Delta^2 w + 464 |\nabla^2 w|^2 - 464 (\Delta w)^2 - 464 \Delta u\ |\nabla u|^2 -464 \langle \nabla w, \nabla |\nabla w|^2 \rangle.
\end{align}
Any solution of (\ref{ELS4}) can be pulled back via stereographic projection to a
solution of (\ref{EulerR4}) with $c = -42$.  Therefore, a corollary of Theorem \ref{Main} is the existence of two distinct rotationally
symmetric solutions of (\ref{EulerR4}) on Euclidean space (with a similar statement for solutions of (\ref{EulerR4tau})). Given this non-uniqueness, it remains an interesting but difficult
problem to classify all entire solutions.  The nonlinear structure of the equations seems to rule out the use of the method of moving planes,
at least in any obvious manner.

In Section \ref{Z} we study rotationally symmetric solutions of (\ref{EulerR4}) with $c = 0$ on $\mathbb{R}^4$ and $\mathbb{R}^4 \setminus \{ 0 \}$, that is, conformal metrics
$g = e^{2w}ds^2$ with $U(g) \equiv 0$.  As in our analysis of the sphere, the problem is reduced to studying the asymptotics of solutions on the cylinder.  We show that there are two families of
periodic solutions, one of which we call {\em Delaunay solutions}, since it includes the cylindrical metric as a limiting case.  The other limiting
case of this family is a solution which we loosely refer to as a {\em Schwarzschild-type} solution.  These solutions are asymptotic to a cone at infinity; see Remark \ref{r:scw} and the example following. We obtain similar results for the half-torsion in Section \ref{ss:gen}.  These examples provides an interesting contrast with our obvious point of comparison, the scalar curvature equation.
\vskip.2in

\subsection{Hyperbolic Space.}  In this paper we study solutions on Euclidean space and the round sphere, but an equally interesting question is the existence of multiple solutions on hyperbolic space.  In \cite{GOR}, the authors show there is an infinite family of rotationally symmetric, complete conformal metrics on the unit ball with constant $Q$-curvature and negative scalar curvature.  In another direction, a \emph{renormalized}
version of the Polyakov formula (\ref{Polyform}) is given in \cite{AAR} for surfaces with cusps or funnels, and the Ricci flow is used to show
the existence of an extremal metric of constant curvature.  It would be very interesting to extend these ideas to four dimensions.
\vskip.2in

\subsection{Organization} The paper is organized as follows:  In Section \ref{s:pf1}  we give the proof of Theorem \ref{Thm1}.  In Section \ref{Z} we consider rotationally symmetric metrics on $\mathbb{R}^4$ with vanishing $U$-curvature.  In Section \ref{Y} we prove the existence of a second critical point on $S^4$ for $F_P$.  In Section \ref{ss:gen} we consider functionals with more general coefficients, and show that the analysis
of Sections \ref{Z} and \ref{Y} apply to the case of the half-torsion.

\section{the proof of Theorem \ref{Thm1}}\label{s:pf1}

The proof of Theorem \ref{Thm1} is elementary, and amounts to gluing in a \emph{bubble} of arbitrary height.  Given $(M^4,g)$, fix a point $p \in M^4$ and let $\{ x^i \}$
denote normal coordinates defined on a geodesic ball $B$ of radius $\rho > 0$ centered at $p$.  Let $\eta \in C_0^{\infty}(M^4)$ be a smooth cut-off function supported
in $B$, and for $\epsilon > 0$ small define
\begin{align*}
w(x) = -\frac{1}{2} \eta \log (\epsilon^2 + |x|^2).
\end{align*}
Using standard formulas for the Laplacian and gradient in normal coordinates, a straightforward calculation gives
\begin{align*}
\int (\Delta w_{\epsilon})^2\ dv &= 4 \omega_3 \log \frac{1}{\epsilon} + O(1), \\
\int (\Delta w_{\epsilon}) |\nabla w_{\epsilon}|^2\ dv &= - 2 \omega_3 \log \frac{1}{\epsilon} + O(1), \\
\int |\nabla w_{\epsilon}|^4\ dv &= \omega_3 \log \frac{1}{\epsilon} + O(1), \\
\log \fint e^{4(w_{\epsilon} - \overline{w}_{\epsilon})}\ dv &= \omega_3 \log \log \frac{1}{\epsilon} + O(1), \\
\int \big [|\Delta w_{\epsilon}| + |\nabla w_{\epsilon}|^2 \big]\ dv &= O(1),
\end{align*}
where $\omega_3$ is the volume of the round $3$-sphere.   Therefore,
\begin{align*}
F_P[w_{\epsilon}] = -24 \omega_3 \log \frac{1}{\epsilon} + O(\log \log \frac{1}{\epsilon}),
& & F_{\tau}[w_{\epsilon}] = -528 \omega_3 \log \frac{1}{\epsilon} + O(\log \log \frac{1}{\epsilon}).
\end{align*}

Letting $\epsilon \rightarrow 0$, we find
\begin{align*}
\inf F_P = -\infty, & &
\inf F_{\tau} = -\infty.
\end{align*}
Replacing $w_{\epsilon}$ with $-w_{\epsilon}$, we also conclude $\sup F_P, \sup F_{\tau} = +\infty$, as claimed.

\section{Metrics of zero $U$-curvature on $\mathbb{R}^4$}\label{Z}

\noindent  In this section we study radially symmetric critical points for the log determinant functional of the Paneitz operator on $\R^4$.  In Section \ref{ss:gen} we will carry out a
similar analysis for the half-torsion.

By (\ref{Idef})--(\ref{IIIdef}) the formula for $\log \det P$ on $\R^4$ is
$$
  L(u) = 18 \int_{\R^4} (\D u)^2 + 64 \int_{\R^4}
  |\n u|^2 \D u + 32 \int_{\R^4} |\n u|^4,
$$
hence we get the following Euler-Lagrange equation:
\begin{equation}\label{eq:eulr40}
    18 \D^2 u + 32 \D \left( |\n u|^2 \right) - 64 div \left(
  \D u \n u \right) - 64 div \left( |\n u|^2 \n u \right) = 0.
\end{equation}
In the space $\mathcal{D}^{2,2}(\R^4)$, the completion of the smooth
compactly supported functions with respect to the Laplace-squared norm,
this functional has a mountain pass structure.

Since we are looking for radial solutions (possibly singular at the origin), it will
be convenient to set up the problem on the cylinder $\mathfrak{C} = \mathbb{R}
\times S^3$ with metric $dt^2 + g_{S^3}$ (conformally equivalent to the flat one).
On $\mathfrak{C}$ one has the identities
$R \equiv 6$ and $Q \equiv 0$, and for $u = u(t)$ we have that
  \begin{equation}\label{eq:Rij}
  R^{ij} \nabla_{ij} u = 0; \qquad \qquad R^{ij} \nabla_i u \nabla_j u = 0.
  \end{equation}
Therefore, the Euler-Lagrange equation becomes
the ODE
\begin{align} \label{eq:eulr4}
9 u^{''''} - 96 u'' (u')^2 + 60 u'' = 0.
\end{align}
Setting $v = u'$ we get
  \begin{equation}\label{eqLode4cons}
  9 v''' - 96 v^2 v' + 60 v' = 0.
  \end{equation}
The latter equation can be integrated, yielding
\begin{equation}\label{eq:vC}
     9 v'' - 32 v^3 + 60 v = \mathcal{C}
\end{equation}
for some $\mathcal{C} \in \R$. This is a Newton equation corresponding to a potential $V_{\mathcal{C}}(v)$ given by
\begin{equation}\label{eq:Vl}
    V_{\mathcal{C}}(v) = - \frac 89 v^4 + \frac{10}{3} v^2 - \frac{\mathcal{C}}{9} v + \frac 23.
\end{equation}
The choice of adding the constant $\frac 23$ in the expression of $V_{\mathcal{C}}$
is for reasons of notational consistency with the next section.

\begin{figure}[h]
\caption{The graph of $V_{\mathcal{C}}$ for $\mathcal{C} = 0, 0 < \mathcal{C}
< 10\sqrt{10}$ and $\mathcal{C} \geq 10\sqrt{10}$}
\begin{center}
 \includegraphics[angle=0,width=6.0cm]{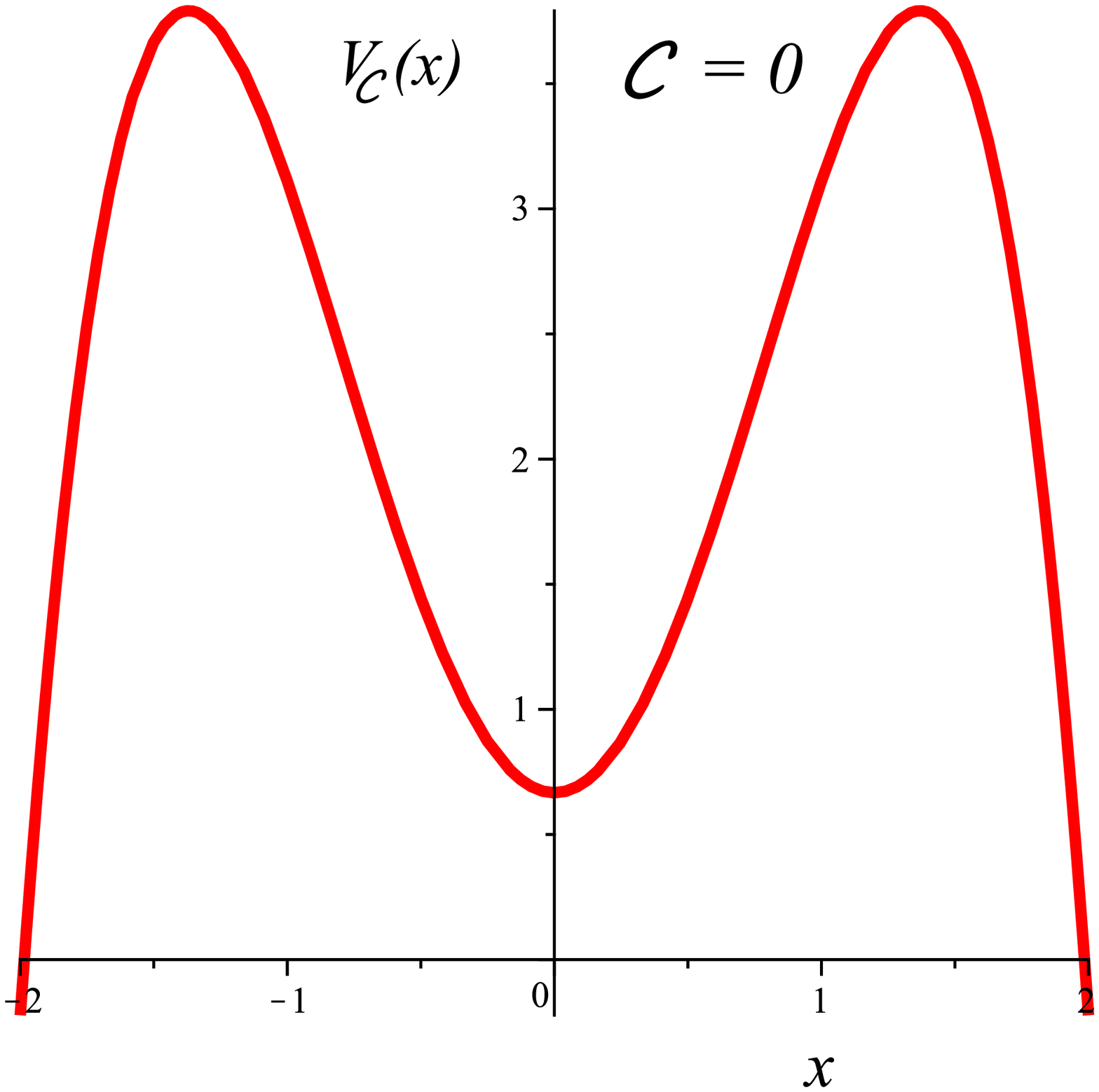} \nonumber \label{fig:0}   \qquad
 \includegraphics[angle=0,width=6.0cm]{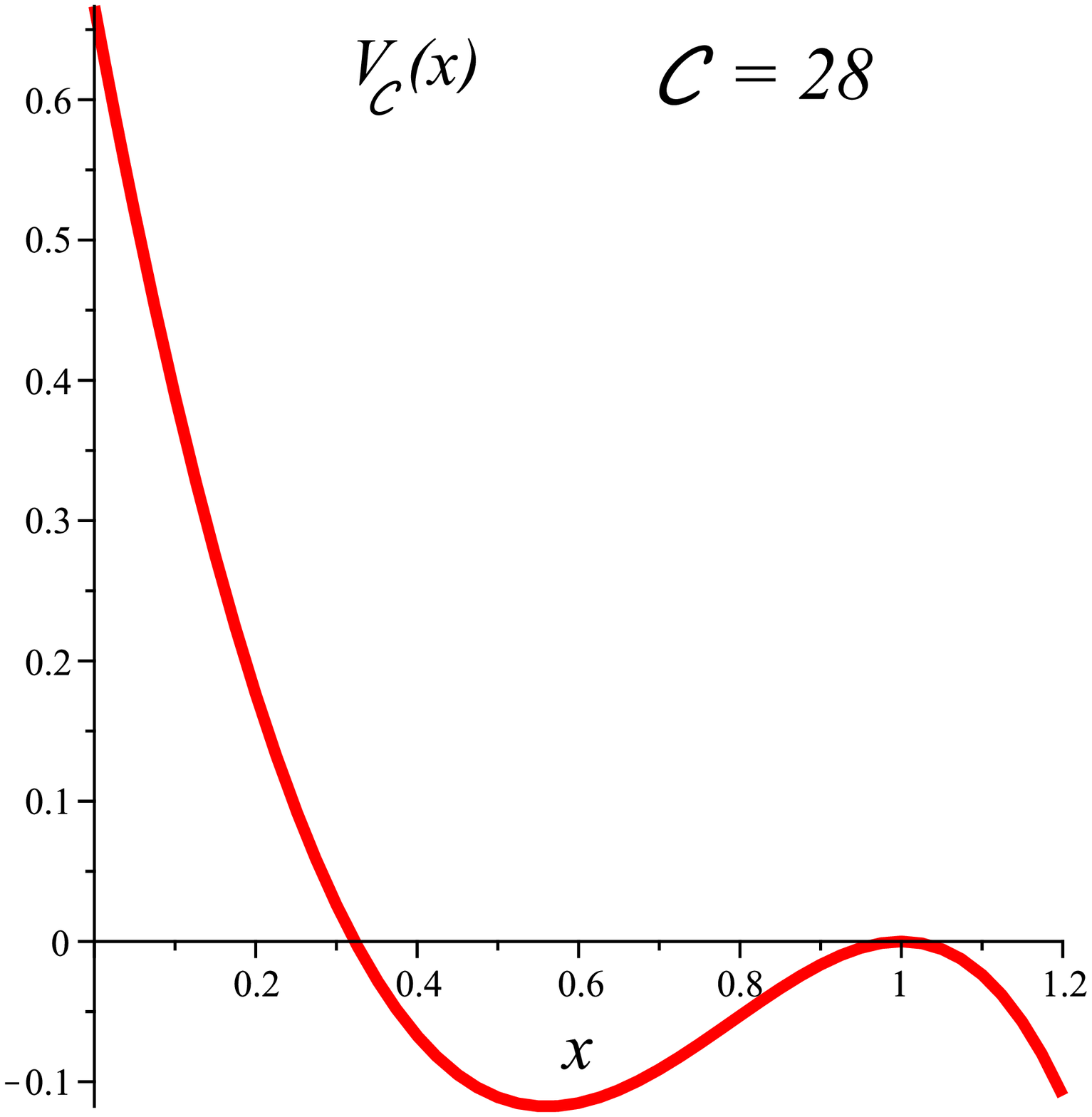} \qquad
 \includegraphics[angle=0,width=6.0cm]{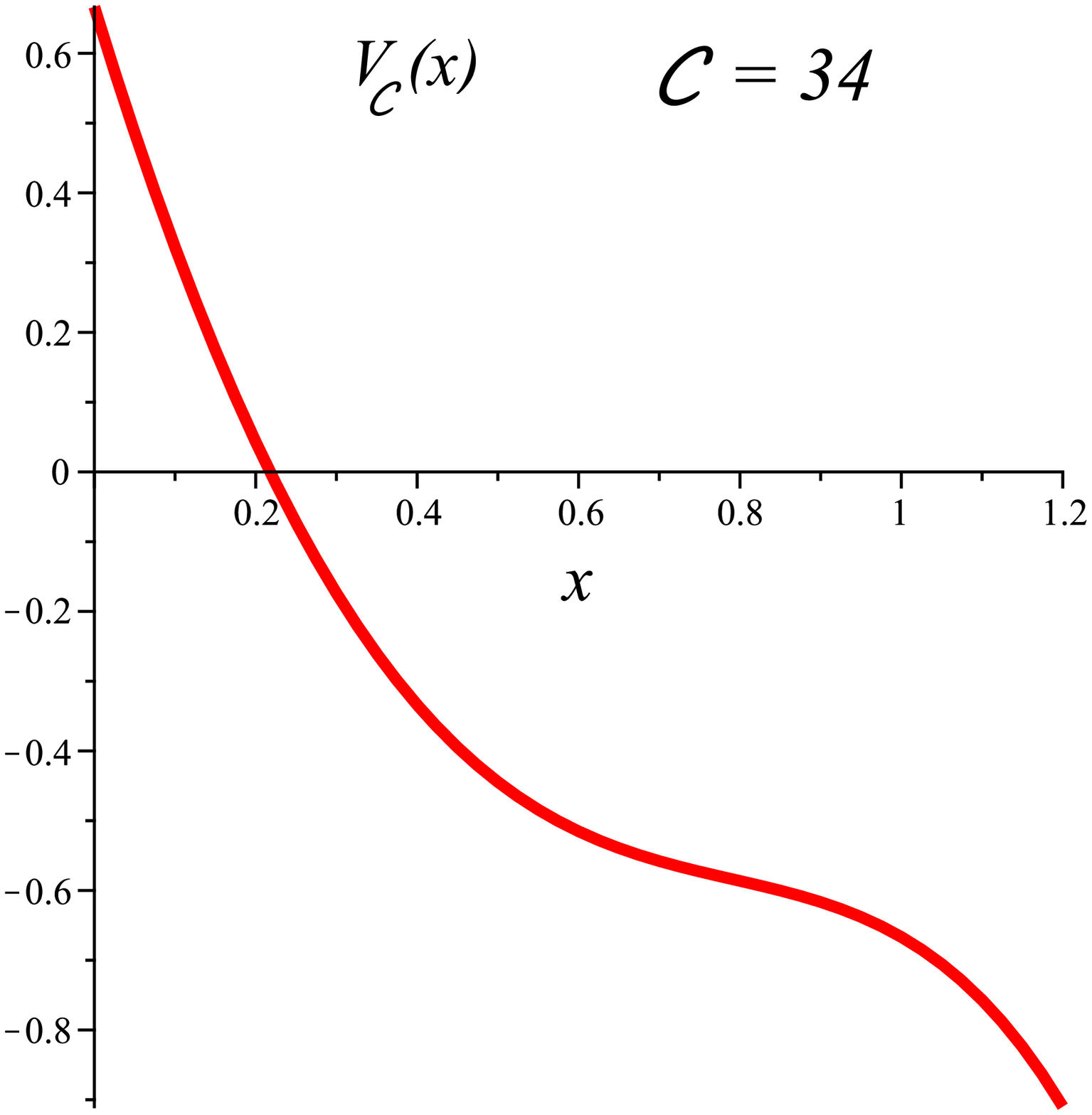} \end{center}
\end{figure}

We divide the analysis into three cases, see Figure \ref{fig:0}. We only consider non
negative values of $\l$, since for $\l < 0$ the situation is symmetric in $v$.
Solutions of \eqref{eq:vC} satisfy the Hamiltonian identity
$$
  \frac 12 (v')^2 + V_{\mathcal{C}}(v) = H,
$$
where $H$ is a constant which depends on the initial data. The latter equation
clearly implies that solutions of \eqref{eq:vC} also satisfy the first order ODE
\begin{equation}\label{eq:1stord}
    v' = \pm \sqrt{2} \sqrt{H - V_{\mathcal{C}}(v)},
\end{equation}
where the $\pm$ sign switches each time $v'$ vanishes and $v'' \neq 0$.

\

\noindent {\bf Case 1: $\mathcal{C} = 0$}

\medskip

\noindent In this situation the potential $V_{\mathcal{C}}$ is even in $v$,
and we can summarize the results in the following proposition.

\begin{proposition}\label{p:delu} (Existence of Delaunay-type solutions) There exist
a one-parameter family of singular solutions $u_\a$ ($\a \in [0,1)$) to \eqref{eq:eulr40},
periodic in $t$, and constants $C_\a > 1$ such that
\begin{equation}\label{eq:uava}
    \frac{1}{C_\a} \frac{(dx)^2}{|x|^2} \leq e^{2 u_\a} (dx)^2 \leq C_\a
  \frac{(dx)^2}{|x|^2} \qquad \quad \hbox{ for all } x \in \R^4.
\end{equation}
For $\a = 0$ we have $e^{2 u_0(x)} \equiv \frac{1}{|x|^2}$.
\end{proposition}

\begin{pf} The existence of a one-parameter family of solutions follows
easily from \eqref{eq:vC}, using the fact that $V_0$ has a reversed double-well
structure with two maxima at $v = \pm \mathrm{v}$, $\mathrm{v} = \frac{\sqrt{30}}{4}$. Their
Hamiltonian energy $H$ ranges in the interval $\left[\frac 23, \frac{91}{24} \right)$.
For $H = \frac 23$ we have a constant solution $v_0 \equiv 0$, corresponding to
the function $u_0$ in the statement of the proposition.

For $H \in \left(\frac 23, \frac{91}{24} \right)$ we obtain a periodic solution
$v_{H}(t)$ oscillating between $- \mathfrak{v}_{H}$ and $\mathfrak{v}_{H}$, where
$0 < \mathfrak{v}_{H} < \mathrm{v}$. Since $v_{H}(t)$ stays uniformly bounded,
we get \eqref{eq:uava} setting $\a = \frac{8}{25} (H - 2/3)$.
\end{pf}

\begin{remark}\label{r:scw} When $H = \frac{91}{24}$ we obtain a heteroclinic
orbit of \eqref{eq:vC} (with $\mathcal{C} = 0$) connecting $- \mathrm{v}$ to $+ \mathrm{v}$.
On $\R^4 \setminus \{0\}$, this corresponds to a solution to \eqref{eq:eulr40} giving
rise to a metric proportional to $r^{- 2 (1 + \mathrm{v})} (dx)^2$ near zero and to
$r^{2 (\mathrm{v}-1)} (dx)^2$ near infinity. These metrics resemble a Schwarzschild
type solution but they are not asymptotically flat near zero or infinity:
asymptotic flatness would correspond to $\mathrm{v}=1$.
\end{remark}

It may help to clarify the preceding remark by considering an explicit example: if we take as
our initial conditions $u'(0) = 0, u''(0) = 5/2$, and $u'''(0) = 0$, then a solution of (\ref{eq:eulr4})
is given by
\begin{align} \label{uex}
u(t) = At + \frac{3}{4} \log \big( 1 + e^{ -\frac{8}{3} A t} \big),
\end{align}
where $A = \sqrt{\frac{15}{8}} > 1$.  Hence,
\begin{align} \label{Uu}
g = e^{2At} \big( 1 + e^{-\frac{8}{3} A t} \big)^{3/2} ( dt^2 + g_{S^3})
\end{align}
is a $U$-flat metric conformal to the cylinder.  Performing the change of variable $r = \frac{1}{A}e^{At}$, we can
write $g$ as
\begin{align} \label{st}
g = [ 1 + O(r^{-8/3})] \big( dr^2 + A^2 r^2 g_{S^3}\big).
\end{align}
Therefore, we see that near infinity, $g$ is asymptotic to a Euclidean cone.
\vskip.1in

\

\noindent {\bf Case 2: $0 < \mathcal{C} < 10\sqrt{10}$}

\

\noindent In this case the potential $V_{\mathcal{C}}$ has two local maxima
$\mathrm{v}_{1,\mathcal{C}} < 0 < \mathrm{v}_{2,\mathcal{C}}$, with
$V(\mathrm{v}_{1,\mathcal{C}}) > V(\mathrm{v}_{2,\mathcal{C}})$.
We have the following proposition.

\begin{proposition}\label{p:2ndr4} For $0 < \mathcal{C} < 10\sqrt{10}$ there exists a
two-parameter family of solutions $u_{\mathcal{C},\a}$ ($\a \in [0,1]$) of \eqref{eq:eulr40}
on $\R^4 \setminus \{0\}$, periodic in $t$, and $C_{\mathcal{C},\a} > 1$,
$\b_{\mathcal{C},\a} \in \R$ such that
\begin{equation}\label{eq:uava2}
    \frac{1}{C_{\mathcal{C},\a}} r^{2 (\b_{\mathcal{C},\a}-1)} (dx)^2
    \leq e^{2 u_{\mathcal{C},\a}} (dx)^2 \leq C_{\mathcal{C},\a}
  r^{2 (\b_{\mathcal{C},\a}-1)} (dx)^2 \qquad \quad \hbox{ for all } x \in \R^4.
\end{equation}
For $\mathcal{C} = 28 \in (0, 10 \sqrt{10})$ and $\a = 1$ we have that $\b_{\mathcal{C},\a} = 1$,
and the metric corresponding to $e^{2 u_{28,1}}$ extends smoothly to (a non flat one on) $\R^4$.
\end{proposition}

\begin{pf}
The proof of the existence part goes exactly as for the previous proposition, with the
difference that when $\a = 1$ we obtain a homoclinic solution (to $\mathrm{v}_{2,\mathcal{C}}$
for $t \to \pm \infty$) instead of a heteroclinic solution.

Let $u_{\mathcal{C},\a}$ be as above and let $\a < 1$: then $v_{\mathcal{C},\a} \equiv
u'_{\mathcal{C},\a}$ oscillates periodically (with period $T_{\mathcal{C},\a}$) between two values
$\tilde{\mathrm{v}}_{\mathcal{C},\a}, \hat{\mathrm{v}}_{\mathcal{C},\a}$, with
$\hat{\mathrm{v}}_{\mathcal{C},\a} > 0$. Suppose that for some $t$ one has
$$
  v_{\mathcal{C},\a}(t) = \tilde{\mathrm{v}}_{\mathcal{C},\a}; \qquad \quad
  v_{\mathcal{C},\a}(t + T_{\mathcal{C},\a}/2) = \hat{\mathrm{v}}_{\mathcal{C},\a}.
$$
Then, from \eqref{eq:1stord} one finds
$$
  T_{\mathcal{C},\a} = 2 \int_{t}^{t + T_{\mathcal{C},\a}/2} ds =
  2 \sqrt{2}
  \int_{\tilde{\mathrm{v}}_{\mathcal{C},\a}}^{\hat{\mathrm{v}}_{\mathcal{C},\a}}
  \frac{dv}{\sqrt{2 (H(\a) - V_{\mathcal{C}}(v))}},
$$
where $H(\a)$ stands for the Hamiltonian energy of the trajectory $v_{\mathcal{C},\a}$.
The number $\b_{\mathcal{C},\a}$ in the statement, which can be taken as the average slope of
$u_{\mathcal{C},\a}$, is given by
$$
  \b_{\mathcal{C},\a} =
  \frac{2}{T_{\mathcal{C},\a}} \int_{t}^{t + T_{\mathcal{C},\a}/2} v_{\mathcal{C},\a}(s)
  ds = \frac{1}{\int_{\tilde{\mathrm{v}}_{\mathcal{C},\a}}^{\hat{\mathrm{v}}_{\mathcal{C},\a}}
  \frac{dv}{\sqrt{2 (H(\a) - V_{\mathcal{C}}(v))}}} \int_{\tilde{\mathrm{v}}_{\mathcal{C},\a}}^{\hat{\mathrm{v}}_{\mathcal{C},\a}}
  \frac{v \, dv}{\sqrt{2 (H(\a) - V_{\mathcal{C}}(v))}}.
$$
For $\a = 1$ then the average of $v$ is simply  $\mathrm{v}_{2,\mathcal{C}}$.

\

When $\mathcal{C} = 28$ one can check that $\mathrm{v}_{2,\mathcal{C}} = 1$, which also
implies $\b_{\mathcal{C},\a} = 1$. For the original solution $u(r)$, this corresponds to
asymptotics of the form $u_{28,1}(r) = C_0 - C_1 r^2 + o(r^2)$ near zero and $u_{28,1}(r)
= C_2 + C_3 r^{-2} + o(r^{-2})$ near infinity. Notice that the flat Euclidean metric
corresponds to $u(t) \equiv t \not\equiv u_{28,1}(t)$. This concludes the proof.
\end{pf}

\

\begin{remark} When $\a$ is small (depending on $\mathcal{C}$) then we can infer
that $\b_{\mathcal{C},\a} > 0$ since $v_{\mathcal{C},\a}$ oscillates near the
local minimum of $V_{\mathcal{C}}$, which is positive.
The same conclusion looks plausible for all $\a \in (0,1]$.
\end{remark}

\

\noindent {\bf Case 3: $\mathcal{C} \geq 10\sqrt{10}$}

\

\noindent In this situation the potential $V_{\mathcal{C}}$ has only one critical
point (a local maximum) $\mathrm{w}_{\mathcal{C}} < 0$ for $\mathcal{C} > 10 \sqrt{10}$, and
two critical points $\mathrm{w}_1 < 0 < \mathrm{w}_2$ for $\mathcal{C} = 10 \sqrt{10}$,
respectively a local maximum and an inflection point. From this structure, one can easily
see that all the globally defined solutions must be constants and coinciding
with some stationary point of $\mathcal{V}_{\mathcal{C}}$.

\section{The proof of Theorem \ref{Main}}\label{Y}

This Section we prove Theorem \ref{Main} for the case of the determinant of the Paneitz operator $F_P$.
In Section \ref{ss:gen} we indicate the necessary changes to prove the result for the half-torsion $F_{\tau}$.

\noindent Recall that the functional determinant for the Paneitz
operator is
$$
F_P[w] = -\frac 14 I - 14 II + \frac 83 III,
$$
whose critical points satisfy the following Euler equation
\begin{align} \label{EL1111} \begin{split}
\mu e^{4w} &=  9 \Delta^2 w +
16 \Delta |\nabla w|^2 - 32 \nabla^i \big[ (\Delta w
+ |\nabla w|^2)\nabla_i w\big]
 - 14 R_{ij}\nabla_i \nabla_j w \\ & \hskip.25in + 10 R\Delta w + 3
\langle \nabla R,\nabla w \rangle - \frac 14 |W|^2 - 14 Q - \frac 83 \Delta R,
\end{split}
\end{align}
where
$$
  \mu = - \frac{\frac 14 \int |W|^2 + 14 \int Q}{\int e^{4w}}.
$$
We will look for solutions on $S^4$ which are radial along some direction and
symmetric with respect to a plane (orthogonal to this given direction), so it will still
be convenient to set up the problem on the cylinder $\mathfrak{C}$, see the beginning of
Section \ref{Z}. Recall that on $\mathfrak{C}$ one has
$R \equiv 6$ and $Q \equiv 0$ and \eqref{eq:Rij}, so if we look for solutions with total
volume equal to $\frac{8}{3} \pi^2$ (the volume one of $S^4$)
from \eqref{EL1111} the Euler-Lagrange equation becomes the ordinary differential equation
\begin{align} \label{ODE4}
9 u^{''''} - 96 u'' (u')^2 + 60 u'' + 42 e^{4u} = 0.
\end{align}
>From the evenness of $u$ we require the initial conditions
\begin{equation}\label{init1}
    \left\{
      \begin{array}{ll}
        u'(0) &= 0, \\
      u'''(0) &= 0.
      \end{array}
    \right.
\end{equation}
Since we need $u$ to lift to a solution on $S^4$ with the correct volume, we also
need the asymptotic conditions
\begin{align} \label{final}
u''(t) \rightarrow 0, \quad u'(t) \rightarrow -1, \quad \int_0^{t} e^{4u}\ ds \rightarrow
\frac{2}{3}, \qquad \ t \rightarrow \infty.
\end{align}

\

\subsection{An auxiliary equation}\label{s:auxeqn}

Using some algebra, we can show that (\ref{ODE4}) reduces to a third
order equation without exponential terms.

\begin{proposition}\label{p:auxeq} Solutions of  \eqref{ODE4} such that
\eqref{init1} and \eqref{final} hold satisfy
\begin{align} \label{init2}
- \frac{9}{2}[u''(0)]^2 + \frac{21}{2} e^{4u(0)}  = 6,
\end{align}
and also the equation
\begin{align} \label{ODE3au}
\frac{9}{4}u'''' - 9 u' u''' - 24 u'' (u')^2 + \frac{9}{2}(u'')^2 +
15 u'' + 24(u')^4 - 30 (u')^2 + 6 = 0.
\end{align}
\end{proposition}

\begin{pf} One can integrate (\ref{ODE4}) and use the initial conditions (\ref{init1}) to get a third order relation:
\begin{align} \label{ODE3}
9 u''' - 32 (u')^3 + 60 u' + 42 \int_0^t e^{4u}\ ds = 0.
\end{align}
Now, multiplying this equation by $u''$ and integrating from $0$ to $t$, integrating by parts in the last term,
and using the initial conditions (\ref{init1}) again, we get
\begin{align} \label{ODE2}
\frac{9}{2} (u'')^2 - \frac{9}{2}[u''(0)]^2 - 8(u')^4 + 30 (u')^2 + 42 u' \int^t_0 e^{4u}\ ds - \frac{21}{2}\big[ e^{4u} - e^{4u(0)} \big] = 0.
\end{align}
Substituting \eqref{final} into (\ref{ODE2}) gives then \eqref{init2}.

Putting this back into (\ref{ODE2}) holds
\begin{align} \label{ODE2a}
\frac{9}{2} (u'')^2  - 8(u')^4 + 30 (u')^2 + 42 u' \int^t_0 e^{4u}\ ds - \frac{21}{2} e^{4u}  + 6 = 0.
\end{align}
Let us now use (\ref{ODE3}) to write
\begin{align*}
42 \int_0^t e^{4u}\ ds = - 9 u''' + 32 (u')^3 - 60 u',
\end{align*}
which implies
\begin{align} \label{1}
42 u' \int_0^t e^{4u}\ ds = - 9 u' u''' + 32 (u')^4 - 60 (u')^2.
\end{align}
Likewise, use the original equation (\ref{ODE4}) to find
\begin{align} \label{2}
- \frac{21}{2} e^{4u} = \frac{9}{4} u^{''''} - 24 u'' (u')^2 + 15 u''.
\end{align}
Substituting these into (\ref{ODE2a}), we eliminate the exponential terms,
and get \eqref{ODE3au}.
\end{pf}

\begin{remark}\label{r:cons}
Putting together \eqref{ODE3} and \eqref{ODE2a} one also finds the conservation law
\begin{align} \label{NT}
9 u''' u' -\frac{9}{2}(u'')^2 - 24 (u')^4 + 30(u')^2 + \frac{21}{2}e^{4u} = 6.
\end{align}
\end{remark}

\

\noindent By Proposition \ref{p:auxeq} and \eqref{init1}, if we let
\begin{align*}
x = x(t) = -u'(t),
\end{align*}
we get the ordinary differential equation
\begin{equation}\label{ODE3a}
    \left\{
  \begin{array}{lll}
    x''' &= & -4 x x'' + \frac{32}{3} x^2 x' + 2 (x')^2
-\frac{20}{3} x' + \frac{32}{3}  x^4  -\frac{40}{3} x^2 + \frac{8}{3},   \\
    x(0) &= & 0,   \\
    x'(0) &=& -u''(0),   \\
    x''(0) &= & 0.
  \end{array}
\right.
\end{equation}
Let us rewrite (\ref{ODE3a}) as a first order system: define
\begin{equation*}
    \left\{
      \begin{array}{lll}
y(t) &= &x'(t), \\
z(t) &= &x''(t).
      \end{array}
    \right.
\end{equation*}
Then (\ref{ODE3a}) is equivalent to
\begin{equation*}
    \left\{
      \begin{array}{lll}
        x' &= &y, \\
y' &=& z, \\
z' &=& -4 x z + \frac{32}{3} x^2 y + 2 y^2
-\frac{20}{3} y + \frac{32}{3}  x^4  -\frac{40}{3} x^2 + \frac{8}{3}.
      \end{array}
    \right.
\end{equation*}
After  some manipulation, we can rewrite this as
\begin{equation}\label{ODE1a}
    \left\{
      \begin{array}{lll}
        x' &= &y, \\
y' &=& z, \\
z' &=& \frac{32}{3}(x-1)\left(x-\frac{1}{2}\right)(x+1)
\left(x+\frac{1}{2}\right) -4 x z + 2 y^2 +
\frac{32}{3} x^2 y  -\frac{20}{3} y,
      \end{array}
    \right.
\end{equation}
with initial conditions
\begin{equation}\label{initxyz}
    \left\{
      \begin{array}{lll}
        x(0) &=& 0, \\
y(0) &=& -u''(0), \\
z(0) &=& 0.
      \end{array}
    \right.
\end{equation}

\

\subsection{Some analysis of \eqref{ODE1a}} One can easily solve for the stationary
points of \eqref{ODE1a}: to begin, putting the first two components equal to zero implies
that $y = z = 0$. Plugging this into the third equation and setting it equal to zero gives
\begin{align*}
\frac{32}{3}(x-1)\left(x-\frac{1}{2}\right)(x+1)\left(x+\frac{1}{2}\right) = 0.
\end{align*}
Therefore,
\begin{align*}
 (x,y,z) \hbox{ is stationary} \quad  \Leftrightarrow \quad (x,y,z) = (\pm 1, 0, 0),
\left(\pm \frac{1}{2}, 0, 0 \right).
\end{align*}
Let $$
p_0 = \left(\frac{1}{2}, 0,0\right), \qquad \qquad  p_1 = (1,0,0),
$$
and let us look at the linearized system at each of these two critical points.
\vskip.2in

\noindent {\bf 1.} At $p_1$, letting $H(x,y,z) = \frac{32}{3}(x-1)
\left(x-\frac{1}{2}\right)(x+1)\left(x+\frac{1}{2}\right) -4 x z + 2 y^2 + \frac{32}{3} x^2 y
-\frac{20}{3} y$, we have
\begin{equation*}
    \left\{
      \begin{array}{ll}
        \frac{\partial H}{\partial x}(p_1) &= 16, \\
\frac{\partial H}{\partial y}(p_1) &= 4, \\
\frac{\partial H}{\partial z}(p_1) &= -4. \\
      \end{array}
    \right.
\end{equation*}
Therefore, the linearized system at $p_1$ is
\begin{equation*}
    \left\{
      \begin{array}{ll}
        x'(t) &= y, \\
y'(t) &= z, \\
z'(t) &= 16x + 4 y - 4z,
      \end{array}
    \right.
\end{equation*}
which we write as
\begin{align*}
\frac{d}{dt}X = {\bf A_1}X,
\end{align*}
with
\begin{equation}\label{eq:A1}
  X
= \left(
    \begin{array}{c}
      x \\
      y \\
      z \\
    \end{array}
  \right); \qquad \qquad \quad
    {\bf A_1} = \left( \begin{array}{ccc}
0 & 1 & 0 \\
0 & 0 & 1 \\
16 & 4 & -4 \end{array} \right).
\end{equation}
The eigenvalues and eigenvectors are
\begin{equation}\label{eq:A1eig}
    {\bf A_1} {\bf v}_1 = 2 {\bf v}_1; \qquad \quad
{\bf A_1} {\bf v}_2 = - 2 {\bf v}_2; \qquad \quad
{\bf A_1} {\bf v}_3 = -4 {\bf v}_3, \qquad \quad
\end{equation}
where
$$
  {\bf v}_1 =
  \left(
      \begin{array}{c}
       1 \\
        2 \\
        4 \\
      \end{array}
    \right) \qquad \quad
  {\bf v}_2 =
   \left(
      \begin{array}{c}
       1 \\
       - 2 \\
        4 \\
      \end{array}
    \right); \quad \qquad {\bf v}_3 =
 \left(
      \begin{array}{c}
       1 \\
        -4 \\
        16 \\
      \end{array}
    \right).
$$
Therefore, $p_1$ is a saddle.
\vskip.2in

\noindent {\bf 2.}   At $p_0$ we have
\begin{equation*}
    \left\{
      \begin{array}{ll}
        \frac{\partial H}{\partial x}(p_0) &= -8, \\
\frac{\partial H}{\partial y}(p_0) &= -4, \\
\frac{\partial H}{\partial z}(p_0) &= -2, \\
      \end{array}
    \right.
\end{equation*}
and the linearized system at this point is
\begin{equation*}
    \left\{
       \begin{array}{ll}
         x'(t) &= y, \\
y'(t) &= z, \\
z'(t) &= -8x  -4y - 2z,
       \end{array}
     \right.
\end{equation*}
which we write as
\begin{align*}
\frac{d}{dt}X = {\bf A_0}X,
\end{align*}
with
\begin{align*}
{\bf A_0} = \left( \begin{array}{ccc}
0 & 1 & 0 \\
0 & 0 & 1 \\
-8 & -4 & -2 \end{array} \right).
\end{align*}
The eigenvalues of this matrix are $\{ -2, -2i , 2i \}$: we will not need the
explicit form of the eigenvectors.

\vskip.2in

\noindent The advantage of looking at system \eqref{ODE3a} instead of the
original equation \eqref{ODE4} is that it is an autonomous one in the
derivatives. Moreover, it includes a one-parameter family of solutions
to \eqref{eqLode4cons}, which is a conservative version of \eqref{ODE4}.

Using our previous notation $(x,y,z)$, \eqref{eqLode4cons} becomes
$K(x,y,z) = 0$, where
$$
  K(x,y,z) = 4 - 6 xz + 3 y^2 + 16 x^4 - 20 x^2.
$$
One can check that the set $\{K = 0\}$ stays invariant for \eqref{ODE3a},
and that solutions on this hypersurface also satisfy \eqref{eqLode4cons}
with $v = - x$. Heuristically, if $u$ attains large negative
values, one might expect that solutions of \eqref{ODE4}-\eqref{final} (and hence
of \eqref{NT}) to behave like those of \eqref{eqLode4cons}. In fact, this is what
we will verify in Subsection \ref{ss:globex} for suitable initial data, see
also Remark \ref{r:approach} below.

We characterize a family of solutions to the first equation of
\eqref{ODE3a} in the following proposition.

\begin{proposition}\label{p:orb} For $\mathcal{C} \in [26,28]$, define
$$
   F_{\mathcal{C}}(x,y,z) = y^2 + 2 V_\mathcal{C}(x); \qquad \qquad
   G_{\mathcal{C}}(x,y,z) = z + \frac{d}{dx} V_\mathcal{C}(x),
$$
where $V_\mathcal{C}$ is given in \eqref{eq:Vl}.
Then for every $\mathcal{C} \in (26,28)$ the system
\begin{equation}\label{eq:systemFG}
    \left\{
      \begin{array}{ll}
        F_{\mathcal{C}}(x,y,z) = 0; & \\
        G_{\mathcal{C}}(x,y,z) = 0 &
      \end{array}
    \right.
\end{equation}
admits a periodic solution $X_{\mathcal{C}}$ which also satisfies
\begin{equation}\label{eq:kkkkkkk}
    x''' =  -4 x x'' + \frac{32}{3} x^2 x' + 2 (x')^2
-\frac{20}{3} x' + \frac{32}{3}  x^4  -\frac{40}{3} x^2 + \frac{8}{3}.
\end{equation}
We get the same conclusions regarding the constant solution $(x(t), y(t), z(t)) \equiv p_0$ when
$\mathcal{C} = 26$, and also for an orbit homoclinic to $p_1$ at $t = \pm \infty$  when $\mathcal{C} = 28$.
\end{proposition}

\begin{pf} Let us first discuss the existence of periodic
solutions of \eqref{eq:systemFG}. The equation
$G_{\mathcal{C}}(x,y,z) = 0$ is a Newton equation for $x(t)$ corresponding to the potential
$V_{\mathcal{C}}$, while the function $F_\mathcal{C}$ stands for (twice) its Hamiltonian energy.
>From the shape of the graph of $V_{\mathcal{C}}$,
see Figure \ref{fig:1}, it is easy to see that periodic solutions with zero Hamiltonian energy
exist for $\mathcal{C} \in (26, 28)$.

\begin{figure}[h]
\caption{The graph of $V_{\mathcal{C}}$ for $\mathcal{C} = 26, 27, 28$}
\begin{center}
 \includegraphics[angle=0,width=6.0cm]{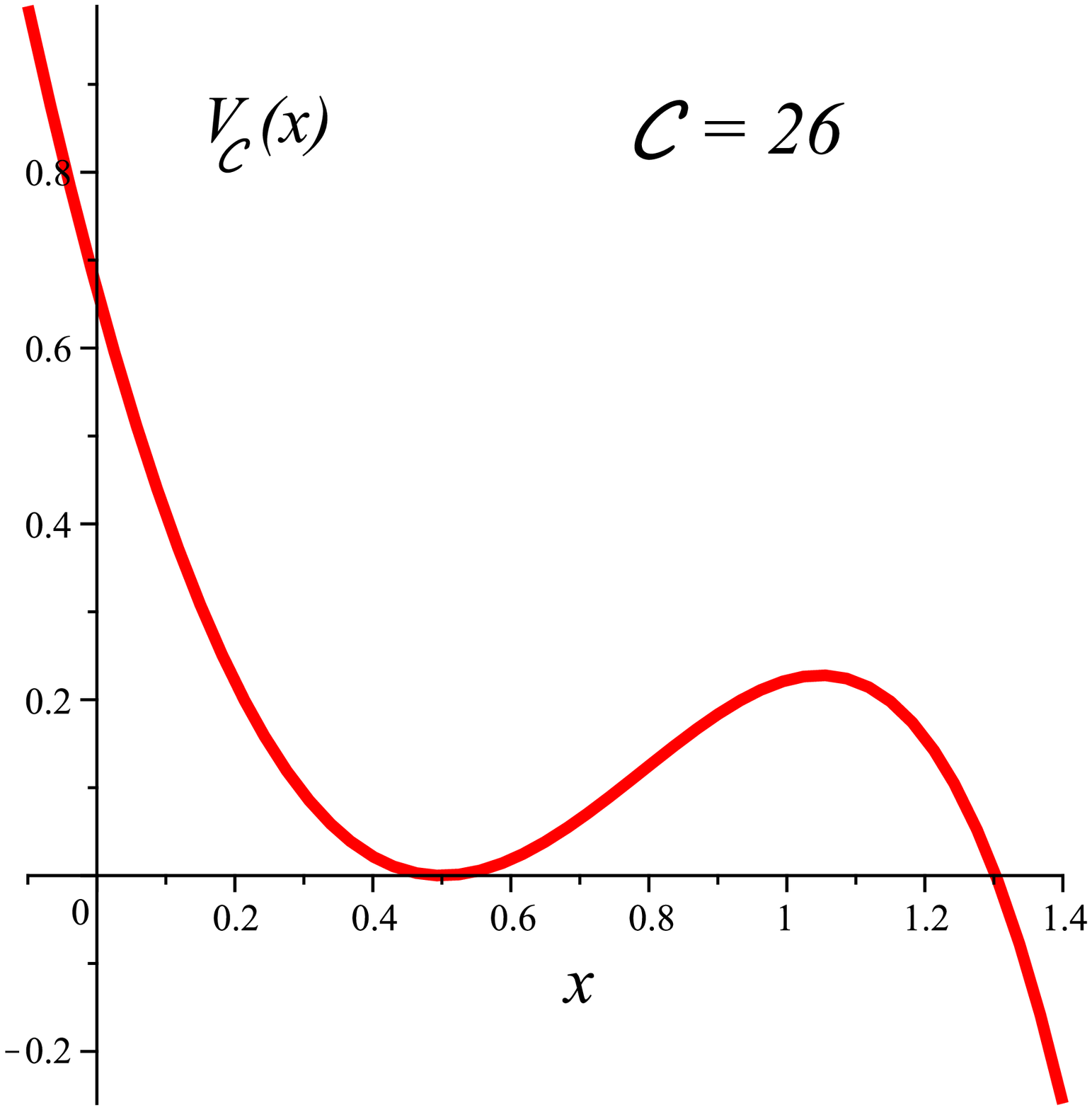} \nonumber \label{fig:1}   \qquad
 \includegraphics[angle=0,width=6.0cm]{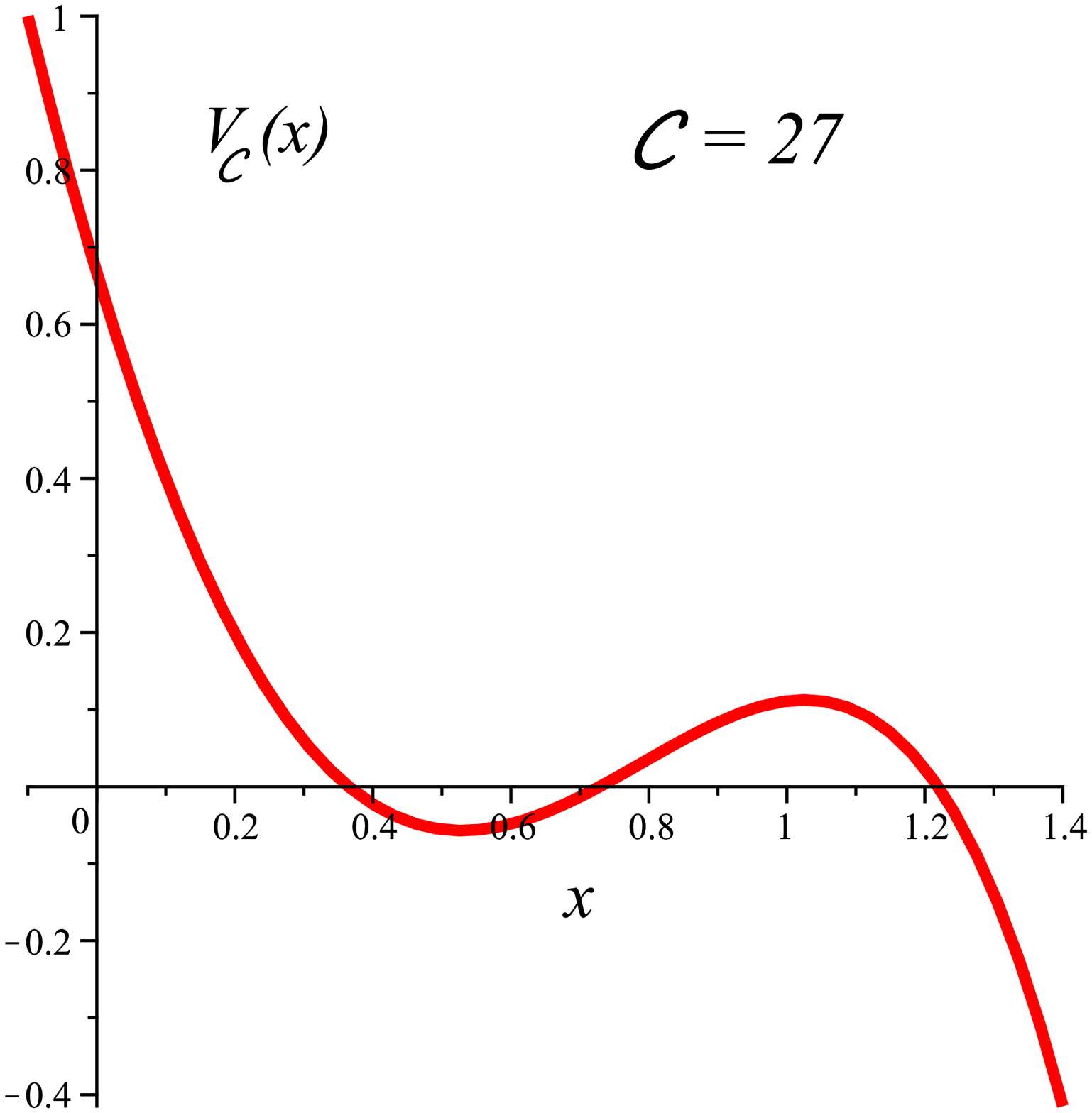} \qquad
 \includegraphics[angle=0,width=6.0cm]{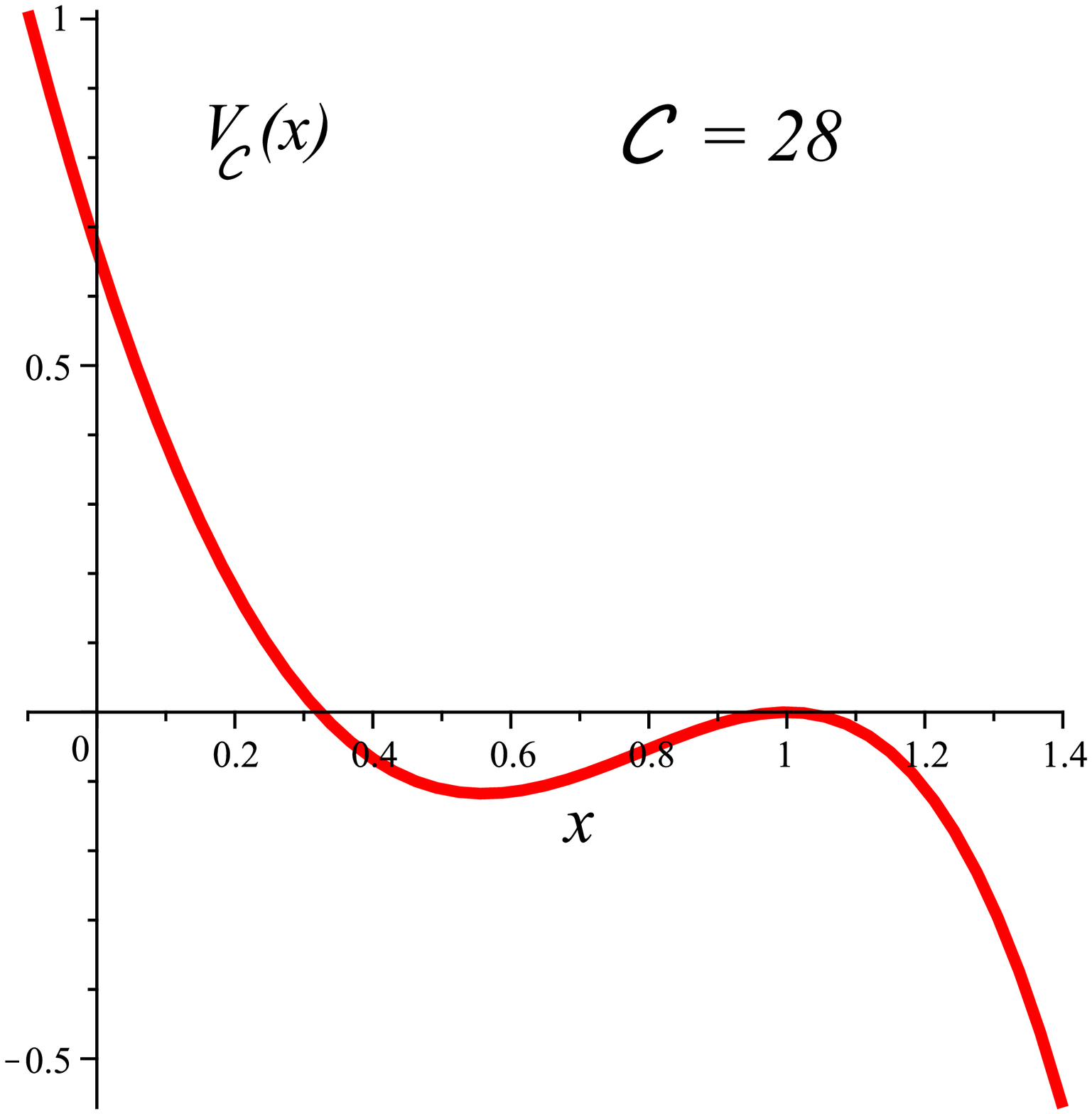} \end{center}
\end{figure}

\noindent The value $\mathcal{C} = 28$
corresponds to a homoclinic solution $X_0 = X_{\mathcal{C} = 28}$ for which
$$
   X_0(t) \to p_1 = (1, 0, 0) \qquad \quad \hbox{ as } t \to \pm \infty.
$$
The value $\mathcal{C} = 26$ instead characterizes the equilibrium point $p_0$ defined above.

By explicit substitution one can easily check that solutions of \eqref{eq:systemFG} also satisfy \eqref{eq:kkkkkkk}.
\end{pf}

\

\begin{figure}[h]
\caption{The topological disc $\mathcal{D}$}
\begin{center}
 \includegraphics[angle=0,width=6.5cm]{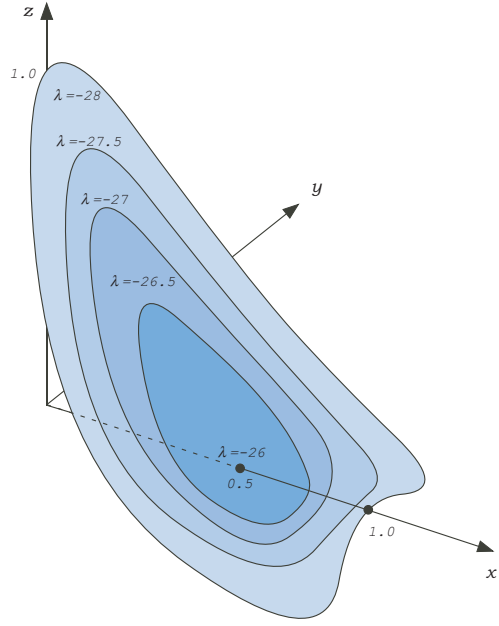} \label{fig:2} \end{center}
\end{figure}

\noindent As $\mathcal{C}$ varies between $26$ and $28$, the trajectories
of $X_{\mathcal{C}}$ foliate a {\em topological disk} $\mathcal{D}$ in $\R^3$ whose boundary
is the homoclinic orbit $X_{\mathcal{C} = 28}$, and whose center is the point
$p_0$, see Figure \ref{fig:2}.

\noindent Let us now go back to equation \eqref{ODE1a}. By some elementary algebra
we obtain the following evolution equations along solutions
\begin{equation}\label{eq:evol}
    \frac{d}{dt} F_\mathcal{C} = 2 y G_\mathcal{C}; \qquad \quad \frac{d}{dt}
    G_\mathcal{C} = \frac 23 K; \qquad \quad \frac{d}{dt} K = - 4 x K.
\end{equation}

\begin{remark}\label{r:approach}
At the points of the disc where $\n K \neq 0$ the last equation in \eqref{eq:evol}
means that the flow is approaching $\mathcal{D}$ (which is contained in the
zero level set of $K$) perpendicularly,
so one could speculate there might exist a positive (in time) invariant set for
\eqref{ODE1a}. This will indeed be proven rigorously in Subsection \ref{ss:globex}.
\end{remark}

\

\noindent We also consider the function
\begin{equation}\label{eq:Q}
  Q(x,y,z) = - 9 z + 32 x^3 - 60 x.
\end{equation}
On the disk $\mathcal{D}$, $Q(x,y,z)$ coincides with
$- \mathcal{C}$, considered as a variable selecting the periodic trajectory. Therefore,
$\mathcal{D} \subseteq \R^3$ can be characterized as
$$
  \mathcal{D} = \{ K = 0 \} \cap \{ -28 \leq Q \leq - 26 \} \cap \{ 0 \leq x \leq 1 \}.
$$
The function $Q$ satisfies the ordinary differential equation
\begin{equation}\label{eq:evolQ}
      \frac{d}{dt} Q = - 6 K.
\end{equation}
Notice that $Q$ coincides with $- 9 G_\mathcal{C} - \mathcal{C}$, so
$\frac{d}{dt} Q$ and $\frac{d}{dt} G_\mathcal{C}$ along
solution have similar expressions.

\

\subsection{Global existence near the spherical metric} \label{ss:globex}
On the cylinder $\mathfrak{C}$, the round metric corresponds to the conformal factor
\begin{align} \label{round}
u_0(t) = - \log \cosh t = \log \left( \frac{2}{e^t + e^{-t}}\right),
\end{align}
which satisfies the initial conditions
\begin{equation}\label{initr}
    \left\{
      \begin{array}{ll}
        u_0'(0) &= 0, \\
u_0''(0) &= -1, \\
u_0'''(0) &= 0.
      \end{array}
    \right.
\end{equation}
The goal of this subsection is to show that for initial data
$$
   \left\{
      \begin{array}{ll}
        x(0) &= 0, \\
   y(0) &= 1-\e, \\
   z(0) &= 0,
      \end{array}
    \right.
$$
with $\e > 0$ small, the solution of \eqref{ODE1a} is globally defined, and
hence also the solution of \eqref{ODE4}.

Let us set
\begin{align} \label{Ndef}
\mathcal{N}[u] = 9 u^{''''} - 96 u'' (u')^2 + 60 u'' + 42 e^{4u},
\end{align}
so that solutions $u$ of (\ref{ODE4}) are characterized by $\mathcal{N}[u] = 0.$
Let $\mathcal{L}$ denote the linearized operator
\begin{align} \label{Lindef}
\mathcal{L}_u\phi = \frac{d}{ds} \mathcal{N}[u + s \phi] \Big|_{s=0}.
\end{align}
If $u = u_0$ is the standard bubble then we simply denote $\mathcal{L}_{u_0}$
by $\mathcal{L}_0$. An easy calculation gives
\begin{align} \label{L0def}
\mathcal{L}_0[\phi] = 9 \phi'''' + [ 60 - 96 (\tanh t)^2
] \phi'' - 192
(\hbox{sech} t)^2 (\tanh t) \phi' + 168 (\hbox{sech} t)^4 \phi.
\end{align}
As $t \rightarrow \infty$, this limits to the equation
\begin{align} \label{L0lim}
\mathcal{L}_0 \phi \sim 9 \phi'''' - 36 \phi'',
\end{align}
so one should expect $\phi$ to be of exponential type at infinity.
\vskip.2in
Indeed, if we linearize the initial conditions on $u$ \eqref{init2}, on
$\phi$ we have to impose
$$
  \phi(0) = - \frac{3}{14}; \qquad \phi'(0) = 0; \qquad \phi''(0) = 1;
  \qquad \phi'''(0) = 0.
$$
An explicit solution is given by the following formula (see Chapter 15 in \cite{AbSteg}
for definitions and properties of hypergeometric functions)
$$
   \phi(t) = -\frac{3}{14} \hbox{hypergeom}\left( \left[\frac 34- \frac{1}{12}
    \sqrt{249},\frac 34
    + \frac{1}{12} \sqrt(249)\right], \frac 12,\cos(2 \arctan(e^t))^2 \right).
$$
By the asymptotics of hypergeometric functions, as $t$ tends to infinity one has
\begin{equation}\label{eq:asyphiA}
      \phi(t) = A (e^{2t} + O(1)); \qquad \qquad \phi'(t) = 2 A (e^{2t} + O(1));
\end{equation}
\begin{equation}\label{eq:asyphiA2}
      \phi''(t) = 4 A (e^{2t} + O(1)); \qquad \qquad \phi'''(t) = 8 A (e^{2t} + O(1))
\end{equation}
for some $A > 0$, where $O(1)$ is a quantity which stays uniformly bounded as $t
\to + \infty$.
\vskip.2in

\noindent We prove first the following result, yielding existence for an interval in the
 variable $t$ which grows as $\e \to 0$, and which relies on a Gronwall type inequality.

\begin{proposition}\label{p:gronwall}
Given $\d > 0$ sufficiently small, there exists $\e_0 > 0$ such that, for $\e \in (0,\e_0)$
the solution of the system with initial data
$$
  (x(0), y(0), z(0)) = (0,1-\e,0)
$$
is defined up to $t = \log \d - \frac 12 \log \e$,  and one has the estimates
\begin{equation}\label{eq:estdelta}
    \left\{
    \begin{array}{ll}
      \left| x(t) - 1 + 2 (e^{-2t} + \e A e^{2t}) \right| \leq \d \e A e^{2t};
     & \\ & \\
    \left| y(t) - 4 (e^{-2t} - \e A e^{2t}) \right| \leq \d \e A e^{2t};  &  \\ &
    \\  \left| z(t) + 8 (e^{-2t} + \e A e^{2t}) \right| \leq \d \e A e^{2t}, &
    \end{array}
  \right. \qquad \hbox{ for } t \in \left[ - \log \d, \log \d - \frac 12 \log \e \right],
\end{equation}
where $A$ is as in \eqref{eq:asyphiA} and \eqref{eq:asyphiA2}.
\end{proposition}

\begin{pf} We can use a Gronwall inequality for the difference between the
true solution and an approximate one. Calling $\var$ the solution to the linearized
equation, we set $x^\e(t) = x^0(t) + \e \var(t) + \tilde{x}^\e(t)$ and then write a
differential inequality for $\tilde{x}^\e$. Recalling that
$$
  X(t) = (x(t), y(t), z(t)),
$$
we write \eqref{ODE1a} in the vector form
$$
  \frac{d}{dt} X(t) = F(X).
$$
We begin by considering the trajectory $X^0$ corresponding to the spherical metric
(for $\e = 0$), which satisfies
\begin{equation}\label{eq:X0}
    \frac{d}{dt} {X}^0 = F(X^0); \qquad \quad X^0(0) = (0,1,0).
\end{equation}
Given a large but fixed $t_0$, by \eqref{round} we have
\begin{equation}\label{eq:asyX^0}
    \left\{
      \begin{array}{ll}
        X^0_1(t_0) & = 1 - 2 e^{-2t_0} + O(e^{-4t_0}); \\
        X^0_2(t_0) & =
  4 e^{-2t_0} + O(e^{-4t_0}); \\
  X^0_3(t_0) & = - 8 e^{-2t_0} + O(e^{-4t_0}).
\end{array}
    \right.
\end{equation}
When we linearize in $\e$ the equation for initial data $X^\e(0) = (0,1-\e,0)$, the
linearized solution satisfies
\begin{equation}\label{eq:dotvar}
    \frac{d}{dt} \var = F'(X^0)[\var],
\end{equation}
with initial conditions $\var(0) = (0,-1,0)$. Recall that, by our previous notation
from Subsection \ref{ss:globex} we have
\begin{equation}\label{eq:varvar}
    \varphi_1 = - \phi'; \qquad \qquad \varphi_2 = -
    \phi''; \qquad \qquad \varphi_3 = - \phi''',
\end{equation}
so from \eqref{eq:asyphiA} and \eqref{eq:asyphiA2} we find
$$
  X^\e_1(t_0) = 1 - 2 e^{-2t_0} - 2 A \e e^{2 t_0} + O(\e^2); \qquad \qquad
   X_2^\e(t_0) = 4 e^{-2t_0} - 4 \e A e^{2t_0} + O(\e^2);
$$
$$
   X_3^\e(t_0) = - 8 e^{-2t_0} - 8 \e A e^{2 t_0} + O(\e^2).
$$
Here, $O(\e^2)$ stands for a quantity bounded by $C_{t_0} \e^2$. We now
choose $\d > 0$ (small but fixed), and then $t_0$ to be the first value of $t$
(depending on $\d$) such that $X^\e_1(t_0) = \frac{\d}{32}$. In this way, we can
write indirectly that $C_{t_0} = C_\d$.

We next set
$$
  X^\e = X^0 + \e \var + \tilde{X}^\e,
$$
and from a Taylor expansion we find
$$
  \left\|F(X^\e) - F(X^0) - \e F'(X_0)[\var] - F'(X^0)[\tilde{X}^\e]
  \var \right\| \leq C_1 \|\e \var + \tilde{X}^\e\|^2,
$$
where $C_1$ is a fixed positive constant (uniformly bounded as long as
the solution lies in a fixed compact set of $\R^3$). Therefore, using
the last formula and some cancelation, we find that
$$
   \left\| \frac{d}{dt} \tilde{X}^\e -  F'(X^0) [\tilde{X}^\e]
   \right\| \leq C_1 \|\e \var + \tilde{X}^\e\|^2.
$$
This implies
$$
  \|\tilde{X}^\e\| \frac{d}{dt} \|\tilde{X}^\e\| = \frac 12 \frac{d}{dt}
  \|\tilde{X}^\e\|^2 = \langle \tilde{X}^\e, \frac{d}{dt} \tilde{X}^\e
   \rangle \leq \langle F'(X^0) [\tilde{X}^\e], \tilde{X}^\e \rangle
  +  C_1 \|\e \var + \tilde{X}^\e\|^2 \|\tilde{X}^\e\|,
$$
and hence
\begin{equation}\label{eq:giulio}
    \frac{d}{dt} \|\tilde{X}^\e\| \leq \frac{\langle F'(X^0) [\tilde{X}^\e],
  \tilde{X}^\e \rangle }{\|\tilde{X}^\e\|} + C_1 \|\e \var + \tilde{X}^\e\|^2.
\end{equation}
Recalling that $F'(0) = {\bf A_1}$, see \eqref{eq:A1}, since $F$ is Lipschitz
by \eqref{eq:asyX^0} one has
$$
   \left\| F'(X^0) - {\bf A_1} \right\| \leq C_1 \|X^0\| \leq \frac{\d}{2}
   \qquad \quad \hbox{ as long as } \quad e^{-2t} \leq \frac{\d}{1000 \, C_1}.
$$
By \eqref{eq:A1eig} this implies
\begin{eqnarray*}
  \frac{\langle F'(X^0) [\tilde{X}^\e], \tilde{X}^\e \rangle
    }{\|\tilde{X}^\e\|} & \leq & \frac{\langle {\bf A_1} [\tilde{X}^\e],
    \tilde{X}^\e \rangle}{\|\tilde{X}^\e\|} + \frac{\d}{2} \|\tilde{X}^\e\| \\
   & \leq &  \left( 2 + \frac{\d}{2} \right) \|\tilde{X}^\e\|
   \qquad  \hbox{ as long as } \quad e^{-2t} \leq \frac{\d}{1000 \, C_1}.
\end{eqnarray*}
>From \eqref{eq:giulio} we then get
$$
  \frac{d}{dt} \|\tilde{X}^\e\| \leq \left( 2 + \frac{\d}{2} \right)
  \|\tilde{X}^\e\|  + C_1 \|\e \var\|^2 + 2 C_1 \|\e \var\| \, \|\tilde{X}^\e\| +
  C_1 \|\tilde{X}^\e\|^2,
$$
which by \eqref{eq:varvar} and \eqref{eq:asyphiA}, \eqref{eq:asyphiA2} yields
$$
  \frac{d}{dt} \|\tilde{X}^\e\| \leq 2 (1 + \d) \|
  \tilde{X}^\e \| + C_1 \|\e \var\|^2 \qquad
\hbox{ as long as } \qquad \left\{
                      \begin{array}{ll}
                        e^{-2t} + \e A e^{2t} \leq \frac{\d}{1000 \,  C_1}; & \\
                        \| \tilde{X}^\e \| \leq \frac{\d}{1000 \,  C_1}. &
                      \end{array}
                    \right.
$$
Therefore, by the asymptotic behavior of $\var$ we have that
\begin{equation}\label{eq:diffineq}
    \frac{d}{dt} \|\tilde{X}^\e\| \leq 2 (1 + \d) \|
  \tilde{X}^\e \| + 128 C_1 \e^2 e^{4t} \qquad
\hbox{ as long as } \qquad \left\{
                      \begin{array}{ll}
                        e^{-2t} + \e A e^{2t} \leq \frac{\d}{1000 \,  C_1}; & \\
                        \| \tilde{X}^\e \| \leq \frac{\d}{1000 \,  C_1}. &
                      \end{array}
                    \right.
\end{equation}
By solving explicitly the associated differential equality, the solution of
\eqref{eq:diffineq} with an initial condition such that
$\|\tilde{X}^\e\|(t_0) \leq C_\d \e^2$ then verifies
$$
\|\tilde{X}^\e\|(t) \leq C_\d \e^2 e^{2(1+\d)(t-t_0)} - \frac{C_1\e^2}{2(1-\d)}
e^{2(1+\d)(t-t_0)} + \frac{C_1 \e^2}{2(1-\d)} e^{4(t-t_0)}
$$
recall, as long as
\begin{equation}\label{eq:constr}
    e^{-2t} + \e A e^{2t} \leq \frac{\d}{1000 \,  C_1} \qquad \quad \hbox{ and }
  \qquad \quad \|\tilde{X}^\e\| \leq \frac{\d}{1000 \,  C_1}.
\end{equation}
We next check the latter condition for $t \in \left[ - \log \d, \log \d - \frac 12
\log \e \right]$. In fact, for $t$ in this range we have that
$$
  e^{-2t} \leq \d^2 < \frac{\d}{1000 \,  C_1}; \qquad \qquad
  \e A e^{2t} \leq A \d^2 < \frac{\d}{1000 \,  C_1},
$$
provided we choose initially $\d$ sufficiently small.

Concerning the second inequality in \eqref{eq:constr}, for
$t \in \left[ - \log \d, \log \d - \frac 12 \log \e \right]$ we get
$$
  \|\tilde{X}^\e\|(t) \leq C_\d \e^2  \left( \frac{\d^2}{\e} \right)^{1+\d}
  + C_1 \e^2 \frac{\d^4}{\e^2} = C_\d \d^{2(1+\d)} \e^{1-\d} + C_1 \d^4 <
  \frac{\d}{1000 \,  C_1},
$$
provided $\d$ is small enough, and if $\e \to 0$. The last estimate also shows
that, for $t$ in the interval $\left[ - \log \d, \log \d - \frac 12 \log \e \right]$
$$
  \|\tilde{X}^\e\|(t) \leq \left( C_\d \d^{2\d-1} \e^{1-\d} + C_1 \d^2 \right)
  \e e^{2t} < \d \e A e^{2t}
$$
for $\d$ sufficiently small, which is the desired conclusion.
\end{pf}

\

\noindent We will show next that, for
suitable initial data close to the ones of the standard bubble,
there exists a globally defined trajectory.

\begin{proposition}\label{p:globexplusasy} For $\e > 0$ small enough the solution $X^\e$ of
\eqref{ODE1a} with initial data
\begin{equation}\label{eq:xeini}
X^\e(0) = (0,1-\e,0)
\end{equation}
is globally defined and there exists $\L_\e \in (-28, -26]$ such that
$$
  K(X(t)) \to 0, \qquad Q(X(t)) \to \L_\e \qquad
  \qquad \hbox{ as } \qquad t \to + \infty.
$$
Moreover, as $t \to + \infty$, $X^\e$ becomes asymptotically periodic.
\end{proposition}

\begin{pf} By Proposition \ref{p:gronwall}, there is $\d$ is sufficiently small
such that, if
 $\e \to 0$, the solution $X^\e$ is defined at least up to  $t = \log \d - \frac{1}{2} \log \e$.
Evaluating it for $\tilde{t}_\e := - \frac 14 \log \e - \frac 14 \log A$ (which is
in the interval where \eqref{eq:estdelta} holds), one has that
$$
  x(\tilde{t}_\e) = 1 - 4 e^{-2\tilde{t}_\e} + 2 e^{-4\tilde{t}_\e} + R_1;
  \qquad \quad y(\tilde{t}_\e) = 4e^{-2\tilde{t}_\e} - 8e^{-4\tilde{t}_\e}
  + R_2;
$$
$$
  z(\tilde{t}_\e) = - 16 e^{-2\tilde{t}_\e} + 32 e^{-4\tilde{t}_\e} + R_3,
$$
where
$$
  |R_i| \leq \d \e A e^{2\tilde{t}_\e} = \d e^{-2 \tilde{t}_\e}.
$$
>From some elementary expansions one finds that
$$
  (Q + 28)(\tilde{t}_\e) = 1320 e^{-4\tilde{t}_\e} + 36 R_1 - 9 R_3 + \tilde{R}_Q;
$$
$$
 K(\tilde{t}_\e) = 688 e^{-4\tilde{t}_\e} + 24 R_1 - 6 R_3 + \tilde{R}_K,
$$
where $|\tilde{R}_Q| + |\tilde{R}_K| \leq C_2 \d \e A$, for a fixed constant $C_2 > 0$. In particular, setting
$$
f(t) = \frac{(Q+28)(t)}{K(t)},
$$
one has $f(t_\e) > \frac 32 + \tilde{C}_\d e^{-2\tilde{t}_\e}$, where $\tilde{C}_\d
\to + \infty$ as $\d \to 0$. Using \eqref{eq:evol} and \eqref{eq:evolQ} one finds that
\begin{equation}\label{eq:odef}
    \frac{d}{dt} f(t) = 4 x(t) f(t) - 6.
\end{equation}
We now estimate the solution from below: setting
$$
  h(t) = \frac 32 + \left( f(\tilde{t}_\e) - \frac 32 \right) e^{3(t-\tilde{t}_\e)},
$$
we show that $h(t)$ is a subsolution of the equation. Since for $t \geq
\tilde{t}_\e$ one has $\e A e^{2t} \geq e^{-2t}$, from the estimates we have on $x(t)$
this would be satisfied if
$$
  3 \left( f(\tilde{t}_\e) - \frac 32 \right) e^{3(t-\tilde{t}_\e)} \leq
  4 (1 - 8 \e A e^{2t}) \left( \frac 32 + \left( f(\tilde{t}_\e) - \frac 32
 \right) e^{3(t-\tilde{t}_\e)} \right) - 6,
$$
namely if
$$
  32 A \e e^{2t} \left( 1 + C_\d e^{3t - 5 \tilde{t}_\e} \right) \leq
  C_\d e^{3t - 5 \tilde{t}_\e} \qquad \hbox{ for } t \geq \tilde{t}_\e.
$$
We claim that this is true for $\tilde{t}_\e \leq t \leq \frac 12 \log \frac{\d^4}{\e A}$.
Notice that this number is smaller than $\log \d - \frac 12 \log \e$, and the estimates
of Proposition \ref{p:gronwall} hold true. We prove that separately
$$
  32 A \e e^{2t} \leq \frac 12 C_\d e^{3t - 5 \tilde{t}_\e} \qquad
  \hbox{ and } \qquad 32 A \e  C_\d e^{5t - 5 \tilde{t}_\e}  \leq
  \frac 12 C_\d e^{3t - 5 \tilde{t}_\e}; \qquad \quad t \leq \frac 12 \log \frac{\d^4}{\e A}.
$$
Taking into account that $e^{- 4 \tilde{t}_\e} = \e A$, the first inequality is equivalent to
$$
  C_\d e^{t - \tilde{t}_\e} \geq 64,
$$
which is true for $\d$ small (recall that $C_\d \to +  \infty$ as $\d \to 0$). The second inequality is  instead equivalent to
$$
  64 A \e e^{2t} \leq 1,
$$
but since $t \leq \frac 12 \log \frac{\d^4}{\e A}$ we have
$$
  64 A \e e^{2t} \leq 64 \d^4,
$$
which is true for $\d$ small. Therefore, we proved  that $h(t)$ is a subsolution.

Hence by comparison we find
$$
  f \left( \frac 12 \log \frac{\d^4}{\e A} \right) \geq h \left( \frac 12
 \log \frac{\d^4}{\e A} \right) \geq \frac 32 + C_\d e^{-5 \tilde{t}_\e}
  \frac{\d^6}{(\e A)^{\frac 32}}.
$$
Using the choice of $\tilde{t}_\e$ then we obtain
$$
  f \left( \frac 12 \log \frac{\d^4}{\e A} \right) \geq \frac 32 + C_\d
  (\e A)^{\frac 54} \frac{\d^6}{(\e A)^{\frac 32}} \geq C_\d \d^6 (\e A)^{- \frac 14}.
$$
This means that
\begin{equation}\label{eq:M}
    \hbox{ for any } M > 0 \hbox{ there exists } t_{\e,M} > \tilde{t}_\e
    \hbox{ such that } f(t) = M.
\end{equation}

\

\noindent Now we notice that
\begin{equation}\label{eq:W}
    W(x,y) := 4 + 20 x^2 - \frac{16}{3} x^4 + 3y^2 - \frac{56}{3} x = K - \frac 23 x (Q + 28)
    = \frac 23 K \left( \frac 32 - x f \right).
\end{equation}
At $t_{\e,M}$ the right hand side is negative. On the other hand
$\{ W < 0 \} \subseteq \R^2$ has a bounded component $W_0$ contained in
$$
  (x,y) \in \left[ \frac{1}{5}, 1 \right] \times \left[ -\frac 35, \frac 35 \right],
$$
(see Figure \ref{fig:3}) and $(x(t_{\e,M}), y(t_{\e,M})) \in W_0$.

\begin{figure}[h]
\caption{The components of $\{ W < 0 \}$ (in white)}
\begin{center}
 \includegraphics[angle=0,width=7.5cm]{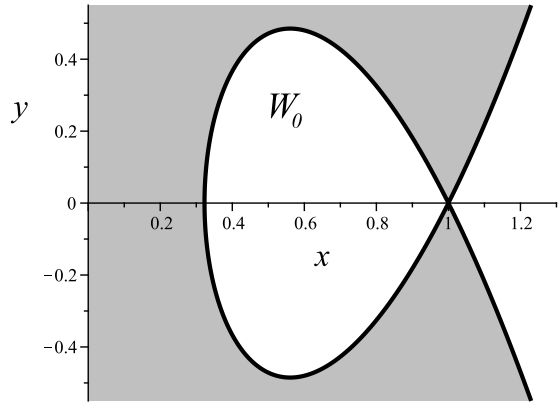} \label{fig:3} \end{center}
\end{figure}

Therefore, as long as $\frac 32 - x(t) f(t) < 0$, we
have a-priori bounds on $x(t)$ and $y(t)$, and $x(t)$ stays positive and
bounded away from zero.
Moreover, $K(t_{\e,M})$ is small positive. Using the expression of
$K$ and the a-priori bounds on $x(t)$ and $y(t)$, one also finds a-priori bounds on $z(t)$,
as long as $\frac 32 - x(t) f < 0$.

If we choose $M > 300$ in \eqref{eq:M}, then $\frac 32 - x(t_{\e,M}) f < 0$ so,
by the bounds on $x(t)$ and by \eqref{eq:odef} $f$ will increase in $t$, so we obtain global existence if $\e > 0$ is sufficiently small.

\

Since $x(t)$ stays bounded, positive and bounded away from zero, by \eqref{eq:evol}
we find immediately that $K(t) \to 0$ as $t \to + \infty$. It remains to prove that
$Q \to \L_\e \in (-28, -26]$ as $t \to + \infty$. Notice that, since $K(0) > 0$ and
since $K(t)$ stays positive, $Q$ is monotone decreasing.

Now, for a small constant $\eta > 0$ and large constant $B > 0$ (to be chosen properly),
we  consider the set
\begin{equation}\label{eq:Oetab}
    \Omega_{\eta, B} := \left\{ 0 \leq K \leq \frac{28 - \eta +Q}{B} \right\}.
\end{equation}
One has that $\n K \neq 0$ on $\{ K = 0 \} \cap \{ Q \in [-28,-26] \}$ (recall that
$\mathcal{D} \subseteq \{ K = 0 \}$ and the third equation in \eqref{eq:evol}), so for
$B$ large and $\eta$ small $\Omega_{\eta, B}$ is a thin neighborhood of the set $\{ K = 0 \}
\cap \{ -28 + \eta \leq Q \leq -26\}$, on the side of $\{K \geq 0\}$.

Using \eqref{eq:evol} and \eqref{eq:evolQ} and the bounds on $x(t)$, one can
check that if $B$ is large then $\Omega_{\eta, B}$ is positive invariant in $t$.
Moreover, from the fact that $f(t) \to + \infty$ as $t \to + \infty$, we can find
$t$ large and $\eta$ small such that $X(t) \in \Omega_{\eta, B}$. Since $Q(t)$ is
monotone decreasing and since $Q \geq - 28 + \eta$ in $\Omega_{\eta, B}$, we obtain that
$Q(t) \to \L_\e \in (-28, 26]$ as $t \to + \infty$, which is the
desired conclusion.  \end{pf}

\

\begin{figure}[h]
\caption{A solution of \eqref{ODE1a} approaching $\mathcal{D}$}
\begin{center}
 \includegraphics[angle=0,width=6.0cm]{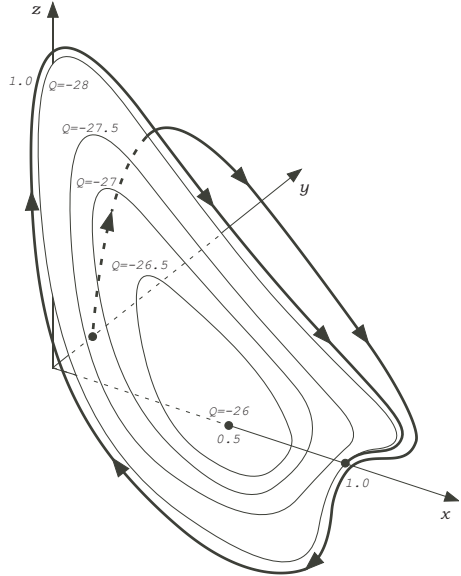}
\label{fig:gg} \end{center}
\end{figure}

\noindent In Figure \ref{fig:gg} a numerical solutions $X^\e$ of \eqref{ODE1a}
is drawn, shadowing one of the periodic orbits in $\mathcal{D}$.

\subsection{A continuity argument}\label{ss:cont} In this subsection we deform the
value of the parameter $\e$ in \eqref{eq:xeini} in order to obtain geometrically
admissible solutions, namely the conditions in \eqref{final}.

\

\noindent Given $\e > 0$, we let $X^\e = (x^\e(t), y^\e(t), z^\e(t))$ denote the
solution of \eqref{ODE1a} with initial condition \eqref{eq:xeini}, and we let $T(\e)$
be the largest number such that $X^\e$ is defined on $[0, T(\e))$. We then let
$\mathcal{E}$ be the family of values of $\e \geq 0$ such that
$$
  \left\{
    \begin{array}{ll}
      T(\e) = + \infty; & \\
      K(X^\e(t)) \to 0 & \hbox{ as } t \to
   + \infty; \\
      Q(X^\e(t)) \to \L_\e \in [-28,-26] & \hbox{ as } t \to
   + \infty.
    \end{array}
  \right.
$$
We also set
\begin{equation}\label{eq:ove}
    \ov{\e} = \sup \left\{ \tilde{\e} \; : \; [0, \tilde{\e}] \subseteq
   \mathcal{E} \right\}.
\end{equation}
First, we show that $\ov{\e}$ is finite.

\begin{lemma}\label{l:bufinite} For $\e > 0$ sufficiently large $T(\e)$ is finite,
and hence $\ov{\e} < + \infty$.
\end{lemma}

\begin{pf} If we define
\begin{equation}\label{eq:mathG}
    \mathcal{G}(t) = x'(t) + 2 x(t)^2 = y(t) + 2 x(t)^2,
\end{equation}
we see that $\mathcal{G}$ satisfies the differential inequality
\begin{equation}\label{eq:diffin}
    \mathcal{G}'' = - \frac{20}{3} \mathcal{G} + \frac{32}{3} x^2 x' + 6 (x')^2
  + \frac{32}{3} x^4 + \frac 83 \geq \frac 83 (\mathcal{G}^2+1) - \frac{20}{3} \mathcal{G},
\end{equation}
and for $t = 0$ we have
\begin{equation}\label{eq:g'zero}
   \mathcal{G}(0) = 1 - \e; \qquad \qquad  \mathcal{G}'(0) = z(0) + 2 x(0) y(0) = 0.
\end{equation}
If we consider the function
$$
  \mathcal{F}(\mathcal{G}, \mathcal{G}') = \frac 12 (\mathcal{G}')^2 - \frac 98
  \mathcal{G}^3 + \frac{10}{3} \mathcal{G}^2 - \frac 83 \mathcal{G},
$$
then by \eqref{eq:diffin} one has that
\begin{equation}\label{eq:ddtf}
    \frac{d}{dt} \mathcal{F}(\mathcal{G}(t), \mathcal{G}'(t)) = \mathcal{G}'(t)
   \left[ \mathcal{G}'' - \frac 83 (\mathcal{G}^2+1) + \frac{20}{3} \mathcal{G} \right].
\end{equation}
For $s > \frac 89$ one can check that the level set $\mathcal{F}(\mathcal{G}, \mathcal{G}') = s$
has only one component, it is symmetric with respect to the $\mathcal{G}$ axis,
it intersects it only once and that
$$
  \left\{ \mathcal{F}(\mathcal{G}, \mathcal{G}') = s \right\} \cap \{ \mathcal{G}'
  \geq 0 \} =  \left\{ (\mathcal{G}, \tilde{F}_s(\mathcal{G})) \; : \; \mathcal{G}
   \in [a_s, + \infty) \right\},
$$
where
$$
  a_s < 0 \hbox{ is decrasing in } s \hbox{ and } a_s \to - \infty \qquad \quad \hbox{ as }
   s \to + \infty;
$$
$$
\tilde{F}_s(\mathcal{G}) > 0,
  \tilde{F}_s(\mathcal{G}) \to + \infty \qquad \quad  \hbox{ as } \mathcal{G} \to + \infty.
$$
Moreover, all these level sets are non degenerate and foliate an open subset
of $\R^2$.

>From this description it follows that if $\e$ is large enough, which implies
that $\mathcal{F}(\mathcal{G}(0), \mathcal{G}'(0))$ is also large,
from \eqref{eq:diffin}, \eqref{eq:g'zero} and \eqref{eq:ddtf} we deduce that
$\mathcal{G}'(t) > \d_\e > 0$ for all $t \in [0, T(\e))$ and that $\mathcal{F}(\mathcal{G}(t), \mathcal{G}'(t))$ increases for all $t \in [0, T(\e))$.

As a consequence, $\mathcal{G}(t)$ becomes large positive with positive
derivative in finite time, so from \eqref{eq:diffin} we deduce that $\mathcal{G}(t)$
must blow up in finite time.
\end{pf}

\begin{lemma}\label{l:G'}
Let $\e \in (0, \ov{\e})$. Then $\mathcal{G}(t)$ and $\mathcal{G}'(t)$ are
uniformly  bounded in $t$.
\end{lemma}

\begin{pf} Similar to the previous proof one can check that for $s < \frac 89$
$$
  \left\{ \mathcal{F}(\mathcal{G}, \mathcal{G}') = s \right\} \cap \{ \mathcal{G}'
  \geq 0 \} \cap \{ \mathcal{G} > 2 \} =  \left\{ (\mathcal{G}, \hat{F}_s(\mathcal{G}))
   \; : \; \mathcal{G} \in [b_s, + \infty) \right\},
$$
where
$$
  b_s > 0 \hbox{ is decrasing in } s \hbox{ and } b_s \to + \infty \qquad \quad  \hbox{ as }
   s \to - \infty;
$$
$$
 \hat{F}_s(\mathcal{G}) > 0,
  \hat{F}_s(\mathcal{G}) \to + \infty \qquad \quad  \hbox{ as } \mathcal{G} \to + \infty.
$$
With the same argument one can prove that if for some $t$ $\mathcal{G}(t)
> 2$ and $\mathcal{G}'(t) > 0$, then there is blow-up in finite time.

\

\noindent As a consequence of this, we deduce that if $\e \in (0, \ov{\e})$ then
$\mathcal{G}$ is uniformly bounded. In fact, since $\mathcal{G}(0)$ is uniformly
bounded for $\e \in (0,\ov{\e})$, if $\mathcal{G}(t)$ becomes large negative for some $t$
by \eqref{eq:diffin} there exists $t_1 > t$ such that $\mathcal{G}(t_1)$ is
large negative and $\mathcal{G}'(t_1) = 0$: we then reason as in the
proof of Lemma \ref{l:bufinite}. If on the other hand $\mathcal{G}(t)$ becomes
large positive for some $t$, then we can argue as before.

\

\noindent  Let us now prove the bounds on $\mathcal{G}'(t)$. If by contradiction
$\mathcal{G}(t)$ stays bounded and $\mathcal{G}'(t)$ becomes large positive, then
$\mathcal{F}(\mathcal{G}(t), \mathcal{G}'(t))$ also becomes large positive, and we can
obtain blow-up in finite time as in the proof of Lemma \ref{l:bufinite}, which
would give a contradiction.

On the other hand, if $\mathcal{G}'(t)$ becomes large negative, it follows from
\eqref{eq:diffin} (arguing as before, but going backwards in $t$) that
$\mathcal{G}(t_2)$ has to be large negative for some $t_2 < t$: we then get a
contradiction from the arguments of the previous paragraph. This concludes the
proof.
\end{pf}

\vskip.2in

\noindent We have next the following result, in which we show that $x^\e(t)$ stays
bounded away from zero for $t$ large enough.

\begin{proposition}\label{p:delta} There exist $T > 0$ and $\d > 0$ such that,
for all $\e \in (0, \ov{\e})$, $x^\e(t) \geq \d$ for $t \geq T$.
\end{proposition}

\

\noindent Since the proof of this proposition is rather long, we begin by stating some
preliminary lemmas, after introducing some useful notation. In the rest of the section
we will always assume that $\e \in (0, \ov{\e})$, and we will often write $X(t), x(t),
\dots$ for $X^\e(t), x^\e(t), \dots$.

\begin{figure}[h]
\caption{The graph of $V_{t}$ for $Q(t) < Q_1, Q(t) = Q_1, Q(t) > Q_1$, $Q_1 \simeq 22.5$}
\begin{center}
 \includegraphics[angle=0,width=6.0cm]{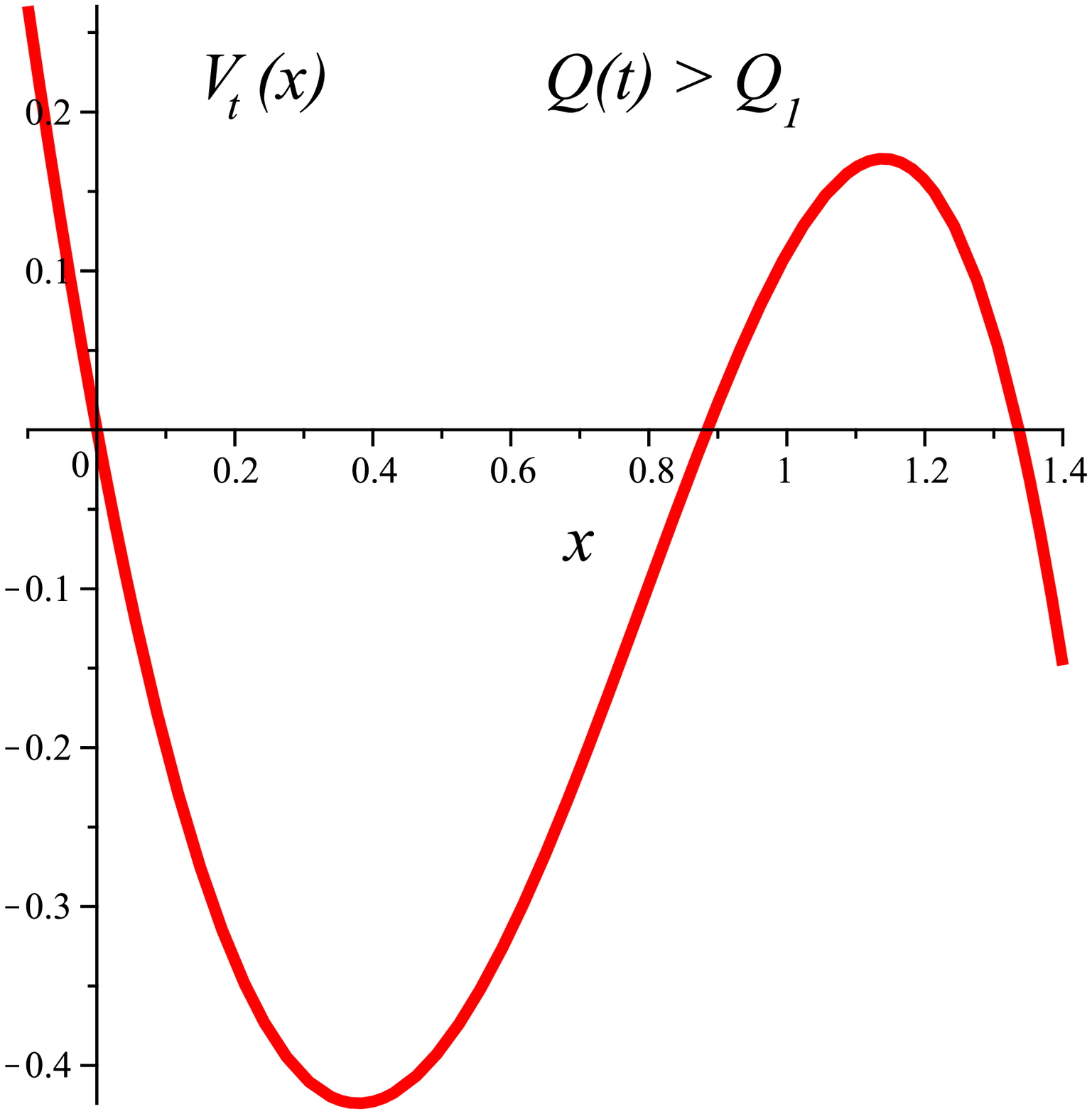} \label{fig:4}   \qquad
 \includegraphics[angle=0,width=6.0cm]{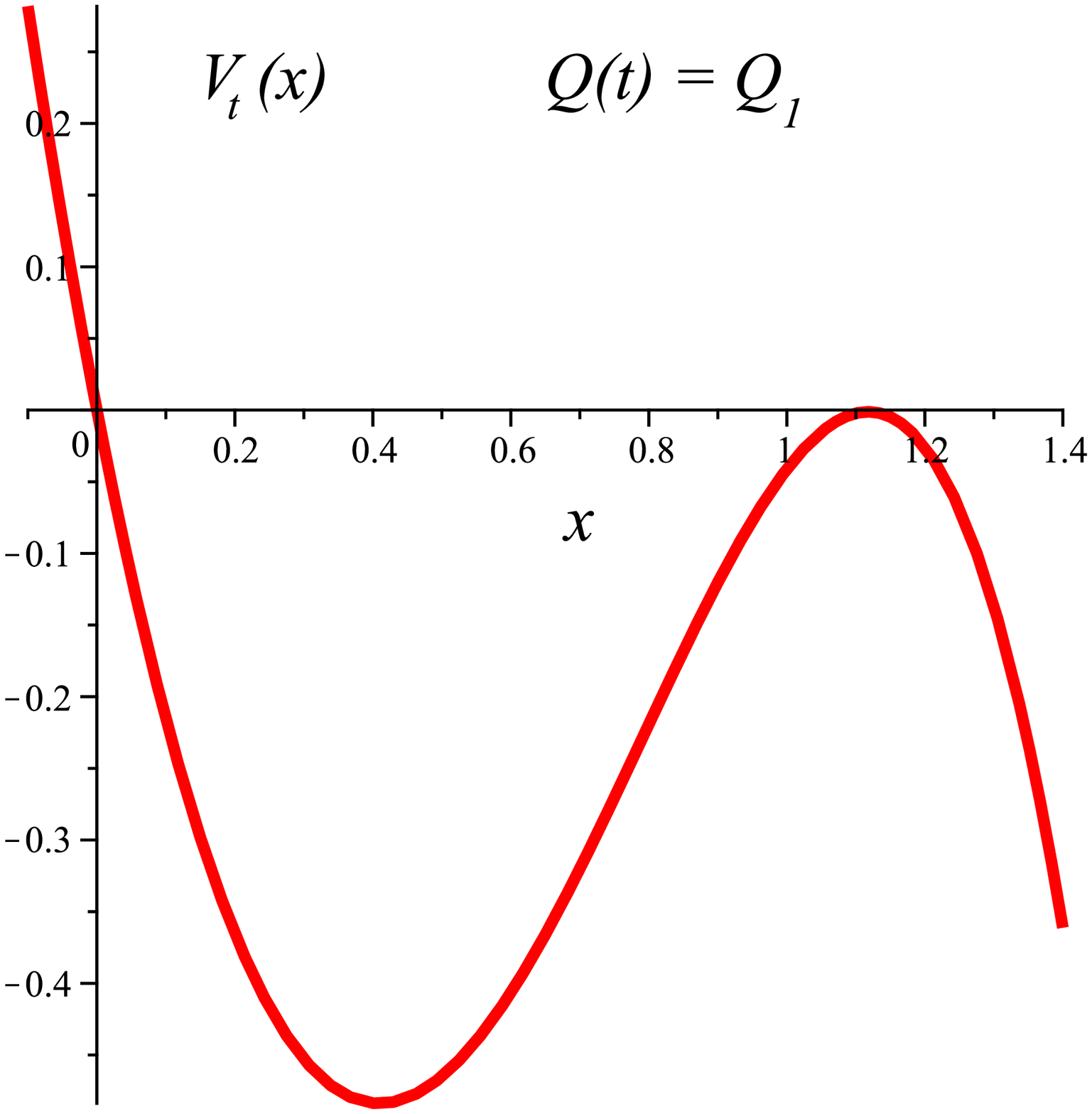} \qquad
 \includegraphics[angle=0,width=6.0cm]{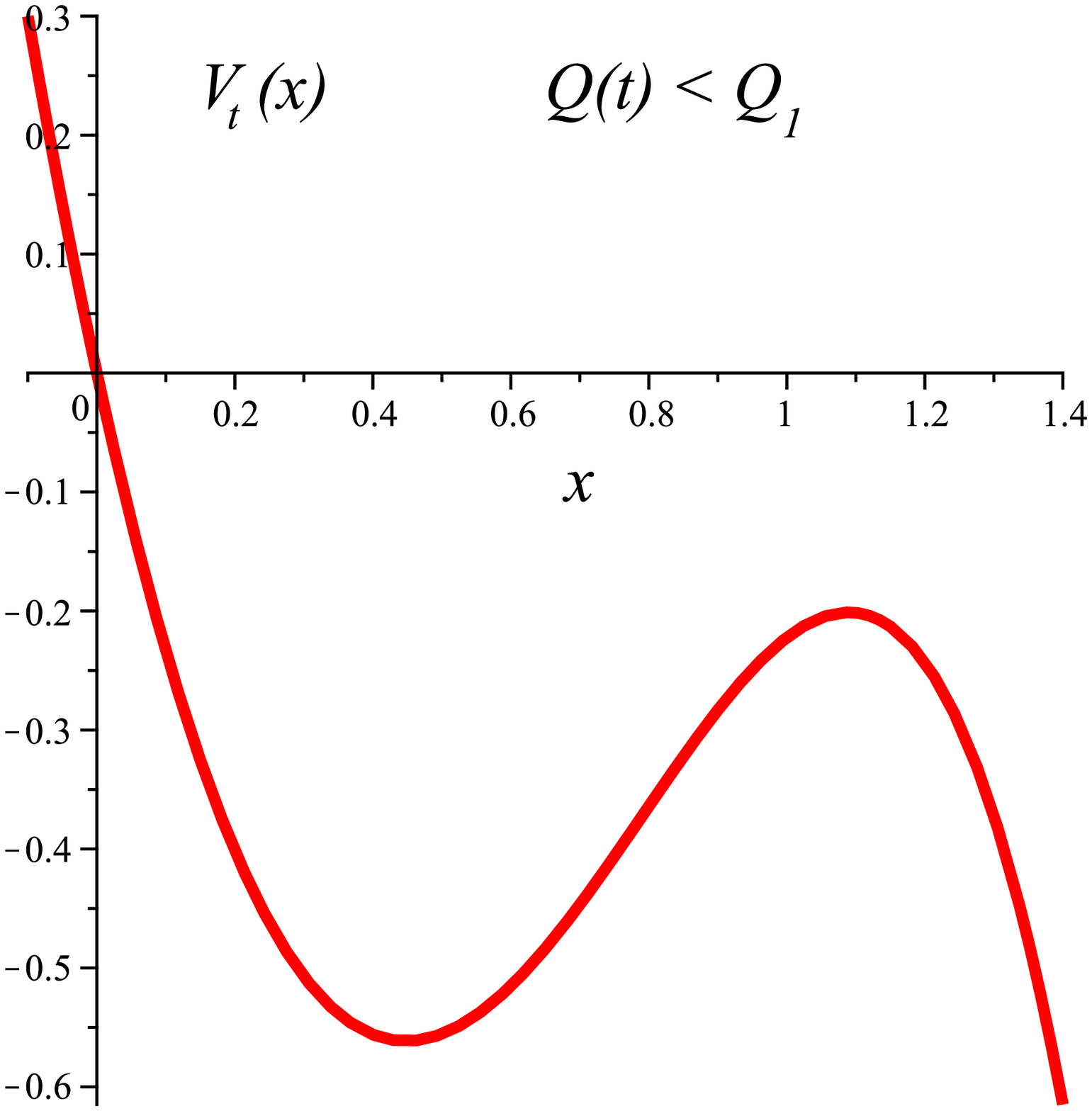} \end{center}
\end{figure}

Recalling the definition of $Q$ in \eqref{eq:Q}, the ordinary differential equation
in \eqref{ODE1a} becomes
\begin{equation}\label{eq:newnew}
    9 x'' = 9 z = - Q(t) - 60 x + 32 x^3.
\end{equation}
This is a Newton equation corresponding to a potential
$V_t$ depending on $t$, which is given by
$$
  V_t(x) = \frac 19 \left( Q(t) x + 30 x^2 - 8 x^4 \right).
$$
For $t = 0$ one has $Q = 0$, so the potential is a reversed double well. Let us examine the
situation for $t > 0$, see Figure \ref{fig:4}.

\begin{lemma}\label{l:xt} For $t > 0$ $V_t(x)$ has negative slope at $x = 0$, and
it has a unique local maximum at $x = \mathfrak{x}_t$ when $x$ positive. Moreover,
we have that
\begin{equation}\label{eq:decrt}
    \frac{d}{dt} \mathfrak{x}_t < 0; \qquad  \frac{d}{dt} V_t(\mathfrak{x}_t) < 0
     \qquad \qquad \hbox{ for all } t,
\end{equation}
and $V_t(\mathfrak{x}_t) < 0$ if $Q(t) \in (-28, -26]$. Furthermore,
for all $t$ either $x(t) < \mathfrak{x}_t$ or $x(t) \geq \mathfrak{x}_t$
and $x'(t) < 0$.
\end{lemma}

\begin{pf} Recalling that $K$ satisfies $\frac{d}{dt} K(t) = - 4 x(t) K(t)$,
see \eqref{eq:evol}, we have that
$$
  K(t) = K(0) e^{- 4 \int_0^t x(s) ds}.
$$
Since $K(0) = 4 + \a^2$, $K(t)$ stays positive for every $t$,
$Q(t)$ decreases to its limit value $\L_\e \in (-28, -26]$
monotonically in $t$. In particular we have that
$$
  \frac{\pa}{\pa x} V_t(x)|_{x=0} < 0 \qquad \quad \hbox{ for } t > 0.
$$
The uniqueness of a local maximum in $x$ for $x > 0$ follows from elementary calculus,
as well as the fact that $V_t(\mathfrak{x}_t) < 0$ if $Q(t) \in (-28, -26]$ and that
$\mathfrak{x}_t \geq 1$ for all $t > 0$. The value $\mathfrak{x}_t$
is defined by the equation
$$
  - Q(t) - 60 \mathfrak{x}_t + 32 \mathfrak{x}_t^3 = 0.
$$
Differentiating with respect to $t$ we obtain
$$
  \left( 96 \mathfrak{x}_t^2 - 60 \right) \frac{d}{dt} {\mathfrak{x}}_t = \frac{d}{dt} Q(t).
$$
The coefficient of $\frac{d}{dt} \mathfrak{x}_t$ in the latter formula is positive by the
fact that $\mathfrak{x}_t \geq 1$, and therefore from $\frac{d}{dt} Q(t) < 0$ and from
some elementary computations we  obtain \eqref{eq:decrt}.

\

\noindent It remains to prove the last statement. Suppose by contradiction that there
exists a first $t_0$ for which
\begin{equation}\label{eq:pospos}
    x(t_0) \geq \mathfrak{x}_{t_0}; \qquad \qquad x'(t_0) \geq 0.
\end{equation}
>From \eqref{eq:newnew}, from the fact that $Q(t)$ is decreasing and from the
fact that $V_t(x)$ has no critical points for $x \geq \mathfrak{x}_{t_0}$,
we deduce that there exists a fixed $\a > 0$ such that
$$
  x''(t) \geq \a (x(t) - \mathfrak{x}_{t_0})^3 \qquad \quad \hbox{ if } x(t)
\geq \mathfrak{x}_{t_0}.
$$
Using the condition \eqref{eq:pospos}  and some comparison arguments
we would then obtain blow-up of $x(t)$ in finite time, which is a
contradiction to the fact that $\e \in (0, \ov{\e})$.
\end{pf}

\

\noindent We next derive some uniform bounds on $X^\e(t)$, together with some useful
consequences.

\begin{lemma}\label{l:bd} There exists a fixed constant $C_0 > 0$ such that $\|X^\e(t)\|
\leq C_0$ for all $\e \in (0, \ov{\e})$ and for all $t > 0$.
\end{lemma}

\begin{pf} We prove first uniform bounds on $x(t)$. We know that $x(0) = 0$, so if
$x(t)$ becomes large positive or large negative the function $\mathcal{G}$ becomes large
positive (either $x(t)$ is large positive and $x'(t) > 0$ for some $t$, or $x(t)$ is
large negative and $x'(t) = 0$ for some $t$: for the latter case, recall that $\e \in (0, \ov{\e})$,
and hence $x(t)$ becomes eventually positive) and has to blow-up in finite time (see
Lemma \ref{l:G'}). This shows uniform bounds on $x(t)$.

Once we have uniform bounds in $x$ we also get uniform bounds in $y = x'$ from those
on $\mathcal{G}$. By Lemma \ref{l:G'} we have that $\mathcal{G}'$
stays uniformly bounded, which implies that also $z = y' = x''$ has to stay uniformly
bounded. Then, using \eqref{ODE1a}, we also get uniform bounds on $x'''$, as required.
This concludes the proof.
\end{pf}

\

\begin{corollary}\label{c:d} There exist $\d_1, \d_2 > 0$ such that
$$
   \frac{d}{dt} Q(t) < - \d_1 < 0 \qquad \quad \hbox{ for all } t \in [0,\d_2] \hbox{
and for all } \e \in (0, \ov{\e}).
$$
Moreover $\frac{d^2}{dt^2} Q(t)$ is uniformly bounded for all $t \geq 0$ and for all
$\e \in (0, \ov{\e})$.
\end{corollary}

\begin{pf} The first statement simply follows from the fact that $\frac{d}{dt} Q(t)
= - 6 K(t)$, that $K(0) = 4 + \a^2$, the continuity of $K(t)$ and from Lemma \ref{l:bd}. The second statement
is immediately deduced from  $\frac{d^2}{dt^2} Q(t) = - 6 K'(t) = 24 x K$ and also from
Lemma \ref{l:bd}. \end{pf}

\

\noindent We next analyze the behavior of solutions when $x(t)$ attains some
small positive value.

\begin{lemma}\label{l:del} There exist $\d_3, \d_4 > 0$ small and $T_0 > 0$, both independent of $\e$,
such that if $x^\e(\ov{t}) = \d_3$ then either
$$
\d_3 < x^\e(s) \leq 3 \qquad \hbox{ for all } s > \ov{t},
$$
or
$$
  \hbox{ there exists } s \in [\ov{t} + \d_4, \ov{t} + T_0] \hbox{ such that } x^\e(s) = \d_3.
$$
\end{lemma}

\begin{pf} First, we show that there exist $T_1 > 0$ large and $\d_3 > 0$ small such that
we have the following implication
\begin{equation}\label{eq:impl}
  \left\{
    \begin{array}{ll}
      t_1 < t_2; &  \\
       x(t_1) = x(t_2) = \d_3; &  \\
      x(t) \geq \d_3 \hbox{ for } t
    \in [t_1, t_2]  &
    \end{array}
  \right.
\qquad \qquad  \Rightarrow \qquad \qquad |t_1 - t_2| \leq T_1.
\end{equation}
In fact, suppose that $x(t) \geq \d_3$ on $[t_1, t_2]$. Since $\frac{d}{dt} K(t) = - 4 x K$,
it means that $K$ shrinks exponentially fast for $t \in [t_1, t_2]$.
Hence, since $K$ is uniformly bounded (and in particular for $t = t_1$) it will
get close to zero if $|t_1 - t_2|$ becomes large. Now notice that
\begin{equation}\label{eq:magicK}
   \frac 16 K = \frac{2}{3} + \frac 12 y^2 + V_t(x)
\end{equation}
and that, by Corollary \ref{c:d}, $V_t$ at $x = 0$ has negative slope (in fact, bounded
away from zero for $t \geq t_1$). Therefore by  $\frac 12 y^2 + V_t(x)
\simeq - \frac 23$, following from \eqref{eq:magicK} and the fact that $K$ is small,
we deduce that $x(t_2)$ cannot approach zero if $K$ is close to zero. This implies
then \eqref{eq:impl} for $\d_3$ small enough.

Let us now prove the statement of the lemma, assuming by contradiction that none of
the two alternatives holds. Let us first suppose that also $x'(\ov{t}) \geq 0$. By \eqref{eq:newnew}
and by Corollary \ref{c:d} we have that $x'(s) > 0$ and $x''(s) > \d_5 > 0$ for $s$ in a right
neighborhood of $\ov{t}$ (of size independent of $\e$).
Therefore, by Lemma \ref{l:bd} (in particular by the bounds on $z$) $x(s) > \d_3$
for $s \in (\ov{t}, \ov{t}+\d_4)$ if $\d_3, \d_4$ are sufficiently small. Since we are disclaiming the
first alternative of the lemma, there will be a first $\tilde{t} > \ov{t}$ for which again
$x(\tilde{t}) = \d_3$. But then we can apply \eqref{eq:impl} to see that we are in the
second alternative.

Suppose now that $x'(\ov{t}) < 0$. By Corollary \ref{c:d} we have that $\frac{d}{dx} V_t(x)$ is negative
and bounded away from zero for $x(t) \leq \d_3$  and for $t > 1$. By \eqref{eq:newnew}, this means that
$x''(s) \geq \d_6 > 0$ for $s > \ov{t}$, as long as $x(s) \leq \d_3$. Therefore (also using the
a-priori bounds in Lemma \ref{l:bd}), we will find $T_2 > 0$ fixed and $\hat{t} \in (\ov{t},
\ov{t} + T_2)$ such that $x(\hat{t}) = \d_3$ and for which $x'(\hat{t}) = \d_7 > 0$, so we end up
in the previous situation ($x'(\ov{t}) > 0$).
\end{pf}

\

\begin{pfn} {\sc of Proposition \ref{p:delta}.} By Lemma \ref{l:del}, (taking $\d = \d_3$ small), the only case we have to exclude is when $x$ equals $\d$ along a sequence $\{t_n\}$, for which $\d \leq t_{n+1} - t_n \leq T_0$. In this case we must have that (see the monotonicity properties in
Lemma \ref{l:xt})
$$
  \lim_{t \to + \infty} V_t(\mathfrak{x}_t) \geq - \d^2
$$
($\d$ is taken small), otherwise by \eqref{eq:newnew} we would deduce blow-up
in finite time (by arguments similar to the proof of Lemma \ref{l:bufinite}).

This means that $Q(t)$ (which is monotone decreasing) stays close to some value $Q_1 > 26$
(again, we are using the smallness of $\d$ and Lemma \ref{l:xt}) on a sequence of intervals $I_n$
of the variable $t$ such that $|I_n| \to + \infty$.
By Corollary \ref{c:d}, we must also have that $\frac{d}{dt} Q(t)$ is small on a sequence of intervals $\tilde{I}_n$
with $|\tilde{I}_n| \to + \infty$, which means (recall the relation $\frac{d}{dt} Q(t) = - 6 K(t)$) that $K$ stays close to zero on the sequence of intervals $\tilde{I}_n$. But this implies that the
function
$$
  \frac 12 y^2 + V_t(x),
$$
the Hamiltonian energy of the trajectory, is negative for $t \in \tilde{I}_n$
(see \eqref{eq:magicK}), so $x(t)$ cannot reach $\d$ for $t \in \tilde{I}_n$.
This concludes the proof.
\end{pfn}

\

\noindent We can now prove the main result of this section.

\begin{proposition}\label{p:adm} The solution $X^{\ov{\e}}$ is globally defined and
satisfies condition \eqref{final}, therefore it is geometrically admissible.
\end{proposition}

\begin{figure}[h]
\caption{A numerical plot of the admissible solution $X^{\ov{\e}}$}
\begin{center}
\includegraphics[angle=0,width=6.0cm]{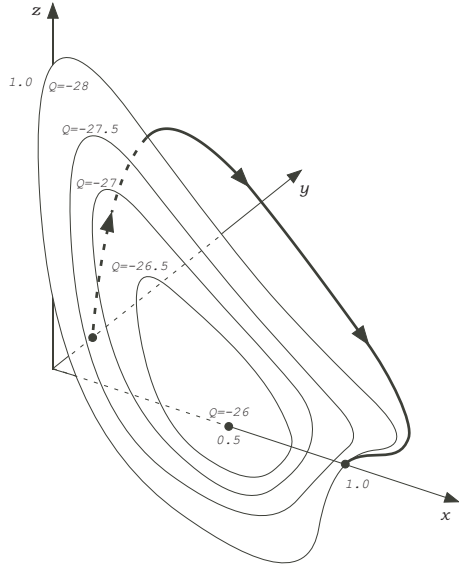} \label{fig:gggg} \end{center}
\end{figure}

\begin{pf} We begin by proving the following claim
\begin{equation}\label{eq:cl1}
    \L_\e \hbox{ is continuous in } \e \hbox{ if } \L_\e \in (-28, -26].
\end{equation}
To see this, let us consider $\e$ such that $\L_\e \in (-28, -26]$, and choose
$\d > 0$ such that $28 + \L_\e > 200 \,
\d > 0$. Let us fix a value of $t$ for which $K(t) < \d$ and $Q(t) - \L_\e < \d$.
>From \eqref{eq:W} and the subsequent arguments one can check that the function $W$ is negative
at $t$, and hence $x(t)$ is positive and bounded away from zero (independently of $\e$ and $\d$
and the times subsequence to $t$).

Choosing now $\tilde{\e}$ for which
$$
  |K^\e(t) - K^{\tilde{\e}}(t)| < \d^2; \qquad \qquad |Q^\e(t) - Q^{\tilde{\e}}(t)| < \d^2,
$$
and using \eqref{eq:evol}, \eqref{eq:evolQ} together with the bounds on $x(t)$ we get
$$
  \frac{d}{dt} \left| K^{\tilde{\e}}(t) - Q^{\tilde{\e}}(t) \right| \leq 200 K^{\tilde{\e}}(t),
$$
which implies, by integration from $t$ to $\infty$, that
$$
  \lim_{s \to + \infty} \left| K^{\tilde{\e}}(s) - Q^{\tilde{\e}}(s) \right| =
   \left| K^{\tilde{\e}}(t) - Q^{\tilde{\e}}(t) \right| + O(\d).
$$
By our choice of $t$ and $\tilde{\e}$ then it follows that
$$
  |\L_\e - \L_{\tilde{\e}}| \leq O(\d),
$$
which implies the continuity of $\L_\e$.

\

\noindent We show next that
\begin{equation}\label{eq:cl2}
   \lim_{\e \nearrow \ov{\e}} \L_{\e} = - 28.
\end{equation}
This follows from the fact that the sets $\Omega_{\eta, B}$ defined in \eqref{eq:Oetab}
are positively invariant in $t$. In fact, suppose that there exist a sequence
$\e_n \nearrow \ov{\e}$ and a fixed $\d > 0$ such that $\L_{\e_n} \geq - 28 + \d$.

By Proposition \ref{p:delta} we deduce uniform (in $\e_n$) exponential decay of $K(t)$
and of $Q(t) - \L_{\e_n}$. This means that we can find $T > 0$ large, $\eta > 0$ small
and $B > 0$ large such that $X^{\e_n}(T) \in \Omega_{\eta/2, 2B}$ for every $n$. But then,
by continuity with respect to the initial data, we an also find $\tilde{\e} > 0$
fixed such that $X^\e \in \Omega_{\eta, B}$ for $|\e - \e_n| \leq \tilde{\e}$.
>From the positive invariance of $\Omega_{\eta, B}$ then we reach a contradiction to
the definition of $\ov{\e}$.

\

\noindent Having \eqref{eq:cl2}, we can now prove the admissibility
conditions \eqref{final}. By Proposition \ref{p:delta} we know that,
if global existence holds,  we have uniform exponential
convergence to one of the periodic orbits (by \eqref{eq:evol} and \eqref{eq:evolQ}
$K$ and $Q$ converge exponentially to their limit values uniformly in $\e \in (0, \ov{\e})$).
When $\L_\e$ approaches $-28$, these periodic orbits have longer and longer period,
and shadow the homoclinic orbit $X_0$ (see the proof of Proposition \ref{p:orb}).
This means that $x(t)$ will be close to $1$ for larger and larger intervals of the
parameter $t$, implying
$$
  \lim_{t \to + \infty} x_{\ov{\a}}(t) = 1.
$$
Therefore our solution is admissible (see Figure \ref{fig:gggg}).
\end{pf}

\section{The case of general coefficients}\label{ss:gen}

\noindent In this section we consider general determinant functionals of the form
\begin{align} \label{FP}
F_A [w] = \gamma_1 I[w] + \gamma_2 II[w] + \gamma_3 III[w].
\end{align}
For convenience we set
\begin{align} \label{alphadef}
\beta  = \frac{\gamma_2}{12 \gamma_3}.
\end{align}
Our goal is to analyze how the arguments in
Sections \ref{Z} and \ref{Y} may be modified as $\beta$ varies (notice that on
locally conformally flat spaces the term $I[w]$ vanishes identically). We are
interested in negative values of $\beta$, since it is for these that $\gamma_2$
and $\gamma_3$ have competing effects. To avoid repetitions, we do not state
explicit results but only limit ourselves to a discussion of the proofs.

\subsection{The zero U-curvature case}\label{ss:uzerogen}

If we study the counterpart of \eqref{eq:eulr40} with a
general choice of the coefficients $\g_i$'s in $\R^4$ and
work on the cylinder $\mathfrak{C} = \R \times S^3$, \eqref{eq:eulr4}
becomes
\begin{equation}\label{eq:euler4gen}
    (1 + \b) u'''' - 6 (u')^2 u'' + (2 - 4 \b) u'' = 0.
\end{equation}
When $\b = -1$ the only solution is $u' = \frac{1}{3}(1-2\b)$, so from now
on we assume that $\b \neq -1$. The case $\mathcal{C} < 0$ is similar
to $\mathcal{C} > 0$, as one can replace $v$ by $-v$.

Integrating \eqref{eq:euler4gen} and setting $v = u'$ we arrive to
$$
  v'' = - V_{\mathcal{C},\b}'(v),
$$
where
$$
   V_{\mathcal{C},\b}(v) = - \frac{1}{2(1+\b)} v^4 + \frac{1-2\b}{1+\b}
   v^2 - \frac{\mathcal{C}}{9} v + \frac 23.
$$
When $\b < -1$, the potential $V_{\mathcal{C},\b}$ is coercive, and periodic
solutions always exist. For $\mathcal{C} = 0$ there are two periodic families
of solutions with $v > 0$ and $v < 0$ respectively, two solutions homoclinic
to zero (giving rise to an asymptotically cylindrical metric), and one family
of periodic changing-sign solutions.

Letting
$$
    \mathcal{C}_\beta = - \frac{12(1-2\b)}{1+\b} \sqrt{\frac{1-2\b}{3}},
$$
a similar qualitative picture, but with a broken symmetry, will persists if
$\mathcal{C} \in (0, \mathcal{C}_\beta)$ (notice that $\mathcal{C}_\beta > 0$
if $\b < -1$). For $\mathcal{C} = \mathcal{C}_\beta$ only one homoclinic
solution will exist, while there will be none for $\mathcal{C} > \mathcal{C}_\beta$.

\

We consider next the case $\b > -1$. When $\mathcal{C} = 0$ we obtain a
one-parameter family of Delaunay type solutions as in Proposition \ref{p:delu}
as well as one heteroclinic solution as in Remark \ref{r:scw}.
When $\mathcal{C} \in (0, - \mathcal{C}_\beta)$ (notice that now $\mathcal{C}_\beta < 0$), the heteroclinic
solution is replaced by a homoclinic solution, while when $\mathcal{C} = - \mathcal{C}_\beta$
only two constant solutions persist.

\subsection{The positive U-curvature case}\label{ss:uposgen}

The Euler equation in this case is given by
\begin{align} \label{ODEgen}
(1 + \beta) u'''' - 6 u'' (u')^2 + (2 - 4 \beta )u'' = c e^{4u},
\end{align}
where the value of $c$ depends on the normalization of $u$.
\vs

Imposing evenness in $t$ and requiring the conditions in \eqref{final} (meaning
the we can lift to a solution on $S^4$) %
%
%
we find
\begin{align} \label{cis}
c = 6 \beta,
\end{align}
so the ODE under interest is
\begin{align} \label{ODEa}
(1 + \beta) u'''' - 6 u'' (u')^2 + (2 - 4 \beta )u'' = 6 \beta e^{4u},
\end{align}
and the integrated version is
\begin{align} \label{ODEb}
(1 + \beta)u''' - 2 (u')^3 + (2 - 4\beta)u' = 6 \beta \int_0^{t} e^{4u}.
\end{align}

\vs

\noindent For the conformal Laplacian $\beta = 1/2$, and the round metric is known
to be the unique even solution. We discuss some features of the values of
$\beta$ smaller than $1/2$, since for $\beta = -7/16$ (corresponding to the
determinant of the Paneitz operator), a second solution exists.

We can now follow the same procedure of reducing the ODE to a third-order system.
The counterpart of \eqref{init2} is
\begin{align} \label{initcon}
-\frac{1}{2}(1 +\beta)[u''(0)]^2 = -\left(2\beta + \frac{1}{2}\right) + \frac{3}{2}\beta e^{4u(0)},
\end{align}
giving the equation
\begin{align} \label{ODEd} \begin{split}
\frac{1}{4}(1 + \beta)u'''' &- (1+\beta)u'''u' + \frac{1}{2}(1+\beta)(u'')^2 - \frac{3}{2}u''(u')^2 \\
& + \frac{1}{4}(2-4\beta)u'' + \frac{3}{2}(u')^4 + (2\beta -1)(u')^2 - \left(2\beta +\frac{1}{2}\right) = 0.
\end{split}
\end{align}

As before, let $x = -u'$, $y = x'$, $z = y'$, we end up with the system
\begin{align} \label{ODEsys} \begin{split}
x' &= y, \\
y' &= z, \\
z' &= -4 x z + \frac{6}{1+\beta} x^2 y + 2 y^2 + 2\left(\frac{2\beta -1}{1 + \beta}\right) y
 \\ & +  \frac{6}{1+\beta}  x^4  + 4\left(\frac{2\beta -1}{1 + \beta}\right) x^2  -
2 \left(\frac{4\beta +1 }{1 + \beta}\right),
\end{split}
\end{align}
with initial conditions
\begin{align} \label{initxyz2} \begin{split}
x(0) &= 0, \\
y(0) &= -u''(0) > 0, \\
z(0) &= 0.
\end{split}
\end{align}

For general $\beta$, we define $K_\beta$ and $Q_\beta$ by
\begin{align} \label{Kdef}
K_\beta = -6xz + 3y^2 + \frac{9}{1+\beta}x^4 - 6\frac{(1-2\beta)}{(1+\beta)} x^2 - 3\frac{(1+4\beta)}{(1+\beta)},
\end{align}
\begin{align} \label{Qdef2}
Q_\beta = -16(1+\beta)z + 32x^3 + 32(2\beta -1)x.
\end{align}
Then along solutions of \eqref{ODEsys} one finds
\begin{align} \label{Kdot}
\frac{d K_\beta}{dt} = -4 x K_\beta,
\end{align}
\begin{align} \label{Qdot}
\frac{d Q_\beta}{dt} = -\frac{32}{3}(1+\beta)K_\beta.
\end{align}

\vs

\noindent Repeating the arguments in the previous subsection one can see that the
limit values of $K_\beta$ and $Q_\beta$ for an admissible solution are $0$ and
$64 \beta$ respectively. The counterpart of \eqref{eq:newnew} is
$$
  x'' = \frac{1}{1+\b} \left[ 2 x^3 + 2 (2\b-1) x - 4 \b \right],
$$
namely a Newton equation with potential
$$
   V_{t,\b}(x) :=  \frac{Q_\b(t)}{16(1+\b)} x -
   \frac{1}{2(1+\b)} x^4 - \frac{2\b-1}{1+\b} x^2.
$$
In the limit $t \to + \infty$, namely when $Q_\b(t)$ tends to $64 \b$,
$V_{t,\b}$ attains a negative maximum at some positive $x$ if and only
if $- 1 < \b < - \frac 14$. For these values of $\b$ then, the above argument
can be repeated with minor changes to get existence of a second solution.
Notice that this applies to the half-torsion case, for which
$\beta = -\frac{31}{58}$.

For $\beta = -1$, $Q_\beta$ has the wrong monotonicity by \eqref{Qdot}, while
$\lim_{t \to + \infty} V_{t,\b}$ has a qualitatively different profile.
For $\b > - \frac 14$ instead, the uniform estimates in Proposition
\ref{p:delta} break down. A numerical simulation indeed indicates that,
although the counterpart of Proposition \ref{p:globexplusasy} holds,
$X^{\ov{\e}}$ is not admissible.

\vskip1in

\bibliography{LogDet_references}

\end{document}